\input ./style/arxiv-general.cfg
\documentclass[aos,MSNbibl,seceqn,dvips]{arximspdf}
\makeatletter
   \@ifpackageloaded{graphicx}{}{\usepackage{graphicx}}
\makeatother
\usepackage{accents}
\usepackage{dcolumn}

\doi{10.1214/15-AOS1322}
\volume{43}
\issue{4}
\pubyear{2015}
\firstpage{1742}
\lastpage{1773}
\docsubty{FLA}

\makeatletter
\newcolumntype{d}[1]{D{.}{.}{#1}}
\newcolumntype{a}[1]{D{;}{\ }{#1}}
\def\sfrac#1#2{#1/#2}

\def\afrac#1#2{#1/(#2)}
\def\vafrac#1#2{(#1)/(#2)}
\def\sklfrac#1#2{(#1/#2)}

\def\sklvafrac#1#2{((#1)/(#2))}
\renewcommand{\mid}{|}
\newcommand{\rrvert}{\vert}
\newcommand{\rrVert}{\Vert}
\newcommand{\llvert}{\vert}
\newcommand{\llVert}{\Vert}
\def\bbR{\mathbb{R}}
\newtheorem{Theorem}{Theorem}[section]
\newtheorem{Lemma}[Theorem]{Lemma}
\newtheorem{Corollary}[Theorem]{Corollary}
\newtheorem{Proposition}[Theorem]{Proposition}
\newproclaim{Assumption}{Assumption}
\newproclaim{Remark}{Remark}[section]
\newproclaim{Example}{Example}[section]

\makeatother

\begin{document}
\begin{frontmatter}

\title{Nonparametric inference in generalized functional~linear models}
\runtitle{Nonparametric inference for functional data}

\begin{aug}
\author[A]{\fnms{Zuofeng}~\snm{Shang}\thanksref{T2}\ead[label=e1]{shang9@purdue.edu}}
\and
\author[A]{\fnms{Guang}~\snm{Cheng}\corref{}\ead[label=e2]{chengg@purdue.edu}}
\runauthor{Z. Shang and G. Cheng}
\affiliation{Purdue University}
\address[A]{Department of Statistics\\
Purdue University\\
250 N. University Street\\
West Lafayette, Indiana 47906\\
USA\\
\printead{e1}\\
\phantom{E-mail:\ }\printead*{e2}}
\end{aug}
\thankstext{T2}{Supported by NSF CAREER Award DMS-11-51692,
DMS-14-18042, Simons Foundation 305266.}

%
\received{\smonth{5} \syear{2014}}
%
\revised{\smonth{2} \syear{2015}}

%
\begin{abstract}
We propose a roughness regularization approach in making nonparametric
inference for generalized functional linear models. In a reproducing
kernel Hilbert space framework, we construct asymptotically valid
confidence intervals for regression mean, prediction intervals for
future response and various statistical procedures for hypothesis
testing. In particular, one procedure for testing global behaviors of
the slope function is adaptive to the smoothness of the slope function
and to the structure of the predictors. As a by-product, a new type of
Wilks phenomenon [\textit{Ann. Math. Stat.} \textbf{9} (1938) 60--62;
\textit{Ann. Statist.} \textbf{29} (2001) 153--193] is discovered when testing the functional linear
models. Despite the generality, our inference procedures are easy to
implement. Numerical examples are provided to demonstrate the empirical
advantages over the competing methods. A collection of technical tools
such as integro-differential equation techniques
[\textit{Trans. Amer. Math. Soc.} (1927) \textbf{29} 755--800;
\textit{Trans. Amer. Math. Soc.} (1928) \textbf{30} 453--471;
\textit{Trans. Amer. Math. Soc.} (1930) \textbf{32} 860--868],
Stein's method
[\textit{Ann. Statist.} \textbf{41} (2013) 2786--2819]
[Stein, \textit{Approximate Computation of Expectations} (1986) IMS]
and functional
Bahadur representation
[\textit{Ann. Statist.} \textbf{41} (2013) 2608--2638]
are employed in this paper.
\end{abstract}

%
\end{frontmatter}

\section{Introduction}

Rapid development in technology makes it possible to
collect measurements intensively over an entire time domain.
This forms the so-called \textit{sample curve}. In functional data analysis,
one may regress the response variable on the sample curve using
(generalized) functional linear models, as in, for example, \cite
{MS05,DPZ12}. Functional principle component analysis (FPCA) is
commonly used for analyzing such models; see, for instance, \cite
{BIW11,YMW05b,HMW06,CH06,HH07,HMV13,L13}. For example, M\"{u}ller and
Stadtm\"{u}ller \cite{MS05} proposed a set of FPCA-based inference
procedures, while Dou et al. \cite{DPZ12} established minimax estimation
rates in a similar framework. The success of these FPCA-based
approaches hinges on the availability of a good estimate of the
functional principal components for the slope function; see \cite
{CY12}. On the other hand, the truncation parameter in the FPCA changes
in a discrete manner, which may yield an imprecise control on the model
complexity, as pointed out in \cite{RS05}. Recently,
Crambes et al. \cite{CKS09}, Yuan and Cai \cite{YC10} and Cai and Yuan \cite{CY12}, among others, have
proposed roughness regularization methods that circumvent the
aforementioned perfect alignment requirement and allow one to
regularize the model complexity in a continuous manner. As far as we
are aware, these works focus mostly on the \textit{estimation or
prediction} in the functional linear models. An exception is the
prediction intervals obtained in \cite{CKS09} under the restrictive
Gaussian errors; see (5.3) therein. However, it is yet unknown how to
handle a broader range of inference problems such as (adaptive)
hypothesis testing for generalized functional linear models in the
above roughness regularization framework.

The major goal of this paper is to systematically conduct asymptotic
inference in the class of generalized functional linear models, which
cover $\ell_2$ regression, logistic regression and exponential family
models. Specifically, we construct confidence intervals for regression
mean, prediction intervals for future response and various statistical
procedures for hypothesis testing. As far as we are aware, all these
inference results are new. In particular, these inference procedures
maintain the modeling and computation flexibility by taking advantage
of the roughness regularization. However, this practical superiority
comes at the price of a much harder theoretical investigation. A key
technical tool we develop in this paper is the \textit{Bahadur
representation for functional data}, which provides a unified treatment
for various inference problems. Due to the involvement of a covariance
operator, we note that this new Bahadur representation is dramatically
different from that recently established in the nonparametric
regression framework \cite{SC13}. In addition, we employ the
integro-differential equation techniques \cite{T1927,TL1928,T1930} to
explicitly characterize the underlying eigen-system that leads to more
transparent inference procedures; see Proposition~\ref{validA3}. As a
side remark, our general theory does not require the Sacks--Ylvisaker
(SY) conditions as in \cite{YC10}, although assuming a pseudo version
of SY conditions (given in Section~S.2) can
facilitate the implementation.

To be more specific, we show that the proposed confidence/prediction
intervals asymptotically achieve
the desirable coverage probability. We also propose a procedure for
testing functional contrast
and show the null limit distribution as a standard normal distribution.
As for testing global behaviors of the slope function, we propose a
penalized likelihood ratio test (PLRT) that achieves the minimax rate
of testing established in \cite{HMV13}. In the particular case of
functional linear models, we observe a new version of the Wilks
phenomenon \cite{W38,FZZ01} arising from PLRT, by which we mean that
the null limit distribution, which is derived as a Chi-square
distribution with diverging degrees of freedom, is free of the true
model parameters. A major advantage of the Wilks type of results is
that we can directly simulate the null limit distribution (without
resorting to bootstrap) in practice. In PLRT, we also point out that
the class of functions in the alternative hypothesis is allowed to be
infinite-dimensional in contrast to the parametric class considered in
\cite{HMV13}. Besides, the rejection region of PLRT is based on the
asymptotic distribution, which makes the procedure more applicable in
general modeling setup, that is, in \textit{generalized} functional linear models.

In reality, the smoothness of the slope function and the structure of
the predictors are typically unknown.
To address this issue, we modify the above PLRT in an ``adaptive''
fashion. Explicitly, we conduct a sequence of standardized PLRT
procedures over multiple smoothness levels, and then use the maximal
one as the new test (after a second standardization). This new testing
method does not rely on prior knowledge of the above two crucial
quantities, and is shown to achieve the minimax rate of testing (up to
logarithm term) established in \cite{HMV13}.
In fact, our adaptive procedures can be viewed as a generalization of
the adaptive Neyman test studied in \cite{F96,FL98} to functional
data. Due to the distinct model structure and test construction, the
Darling--Erd\H{o}s theorem used in \cite{F96,FL98} is no longer
applicable. Instead, we adapt the impressive and powerful Gaussian
approximation tool recently proposed in \cite{CCK13} to show that in
both Gaussian and sub-Gaussian settings, the null limit is a type of
extreme value distribution. Our adaptive testing procedures differ from
the FPCA-based tests such as those considered in \cite{HMV13,L13} in
two ways: (i) our tests work for non-Gaussian models; (ii) our tests
provide an asymptotic null limit distribution, from which the correct
test size can be achieved. Besides, our tests do not require the
``eigen-gap'' condition in the FPCA literature, as in, for example,
\cite{L13}. Simulation results demonstrate the advantages of our
methods in terms of desirable sizes and powers. In particular, we
observe that PLRT is more powerful than the adaptive testing
procedures. This is reasonable since PLRT incorporates prior knowledge
on smoothness of the covariance and reproducing kernels. However, their
difference quickly vanishes when the sample size is large or the signal
strength is strong.
%

The rest of this paper is organized in the following way. In
Section~\ref{secpreliminary}, basic assumptions on model and parameter space
are given. Section~\ref{secbrfd} presents the key technical device of
this paper:
Bahadur representation for functional data. In Section~\ref
{secasympCI}, asymptotically valid confidence intervals
for regression mean and prediction intervals for future response are
constructed.
In Section~\ref{sechypothesistesting}, a procedure for testing
functional contrast and a global testing for the slope function, that
is, PLRT, are established. Theoretical properties are also
demonstrated. Section~\ref{secadaptivePLRT} contains two adaptive
testing procedures for either Gaussian or sub-Gaussian errors. Their
null limit distributions and minimax properties are carefully examined.
A simulation study is provided in Section~\ref{secsimulation}. The
generalized cross validation (GCV) is used to select the roughness
penalty parameter in the simulations. Section~\ref{secdisc} discusses
the technical connection between our work and \cite{SC13}. All
technical proofs are deferred to the Supplementary Material \cite{SCFDA}.

\section{Preliminaries}\label{secpreliminary}
\subsection{Model assumptions}
Suppose the data $(Y_i,X_i(t))$, $i=1,\ldots,n$, are i.i.d. copies of
$(Y,X(t))$,
where $Y$ is a univariate response variable taking values in $\mathcal{Y}$,
a subset of real numbers, and $X(t)$ is a real-valued random predictor
process over $\mathbb{I}=[0,1]$.
Consider the following generalized functional linear model:
%
\begin{equation}
\label{gflm} \mu_0(X)\equiv E\{Y\mid X\}=F \biggl(
\alpha_0+\int_0^1 X(t)\beta
_0(t)\,dt \biggr),
\end{equation}
where $F$ is a known link function,
$\alpha_0$ is a scalar and $\beta_0(\cdot)$ is a real-valued
function. The conditional mean w.r.t. $X=X(\cdot)$ can be understood
as a function of a collection of random variables
$\{X(t)\dvtx  0\le t\le1\}$ throughout the paper. Let $\beta\in
H^m(\mathbb{I})$, the $m$-order Sobolev space
defined by
\begin{eqnarray*}
H^m(\mathbb{I})&=&\bigl\{\beta\dvtx \mathbb I\mapsto\mathbb R\mid
\beta^{(j)}, j=0,\ldots,m-1,
\\
&&\hspace*{4pt}\mbox{are absolutely continuous, and }\beta^{(m)}\in L^2(
\mathbb{I}) \bigr\}.
\end{eqnarray*}
Therefore, the unknown parameter $\theta\equiv(\alpha,\beta)$
belongs to $\mathcal{H}\equiv\bbR^1\times H^m(\mathbb{I})$. We
further assume $m>1/2$ such that $H^m(\mathbb{I})$ is a reproducing
kernel Hilbert space.

In this paper, we consider a general loss function $\ell(y;a)$ defined
over $y\in\mathcal{Y}$ and $a\in\bbR$, which covers two important
classes of statistical models: (i)~$\ell(y; a)=\log p(y; F(a))$, where
$y\mid x\sim p(y;\mu_0(x))$ for a conditional distribution $p$;
(ii)~$\ell
(y; a)=Q(y; F(a))$, where $Q(y; \mu)\equiv\int_{y}^\mu
(y-s)/\mathcal V(s)\,ds$ is a quasi-likelihood with some known
positive-valued function $\mathcal{V}$ satisfying $\mathcal V(\mu
_0(X))=\operatorname{Var}(Y\mid X)$; see \cite{W74}. Note that these two criterion
functions coincide under some choices of~$\mathcal V$. The regularized
estimator is given by
%
\begin{eqnarray}
\label{estimation} \qquad &&(\widehat{\alpha}_{n,\lambda},\widehat{\beta
}_{n,\lambda
})\nonumber
\\
&&\qquad = \arg\sup_{(\alpha,\beta)\in\mathcal{H}}\ell_{n,\lambda
}(\theta)
\\
&&\qquad \equiv \arg\sup_{(\alpha,\beta)\in\mathcal{H}} \Biggl\{ \frac
{1}{n}\sum
_{i=1}^n\ell\biggl(Y_i;
\alpha+\int_0^1 X_i(t)\beta(t)\,dt
\biggr)-(\lambda/2)J(\beta,\beta) \Biggr\},\nonumber
\end{eqnarray}
where $J(\beta,\widetilde{\beta})=\int_0^1\beta^{(m)}(t)\widetilde
{\beta}^{(m)}(t)\,dt$ is a roughness penalty. Here, we use $\lambda/2$
to simplify future expressions. In the special $\ell_2$-regression,
Yuan and Cai \cite{YC10} study the minimax optimal estimation and
prediction by assuming the same roughness penalty.

We next assume the following smoothness and tail conditions on $\ell$.
Denote the first-, second- and third-order derivatives of $\ell(y;a)$
w.r.t. $a$ by $\dot\ell_a(y;a)$, $\ddot{\ell}_a(y;a)$ and $\ell
'''_a(y;a)$, respectively.

%
\begin{Assumption}\label{A1}
\textup{(a)}  $\ell(y;a)$ is three times continuously
differentiable and strictly concave w.r.t $a$.
There exist positive
constants $C_0$ and $C_1$ s.t.,
%
\begin{eqnarray}
\label{A1aeq1} E \Bigl\{\exp\Bigl(\sup_{a\in\bbR}\bigl\llvert
\ddot{\ell}_a(Y;a)\bigr\rrvert/C_0 \Bigr)\big|  X \Bigr\}
&\le& C_1,
\nonumber\\[-8pt]\\[-8pt]\nonumber
E \Bigl\{\exp\Bigl(\sup_{a\in\bbR} \bigl\llvert
\ell'''_a(Y;a)\bigr\rrvert
/C_0 \Bigr) \big|  X \Bigr\}&\le& C_1,\qquad\mbox{a.s.}
\end{eqnarray}

(b) There exists a positive constant $C_2$ s.t.,
\[
C_2^{-1}\le B(X)\equiv-E \biggl\{\ddot{
\ell}_a \biggl(Y;\alpha_0+\int_0^1X(t)
\beta_0(t)\,dt \biggr)\big| X \biggr\}\le C_2\qquad\mbox{a.s.}
\]
In addition, $X$ is weighted-centered in the sense that $E\{B(X)X(t)\}
=0$ for any $t\in\mathbb{I}$.

(c)
$\epsilon\equiv
\dot{\ell}_a (Y;\alpha_0+\int_0^1X(t)\beta_0(t)\,dt )$ satisfies
$E\{\epsilon\mid X\}=0$ and $E\{\epsilon^2\mid X\}= B(X)$, a.s.
\end{Assumption}

The weighted center condition in Assumption~\ref{A1}(b) is only used to simplify our technical arguments. Actually, it
always holds after a simple data transformation; see the Supplementary
Material~\cite{SCFDA}, Section~S.1. Next, we give three
examples to illustrate the validity of Assumption~\ref{A1}.

%
\begin{Example}[(Gaussian model)]\label{exa1pss}
In the functional linear models under Gaussian errors, that is,
$Y=\alpha_0+\int_0^1 X(t)\beta_0(t)\,dt+v$ and $v\mid X\sim N(0,\sigma
^2)$, we can easily verify Assumption~\ref{A1} with $B(X)=\sigma
^{-2}$ and $\epsilon=v/\sigma^2$ given that $E\{X(t)\}=0$.
\end{Example}

%
\begin{Example}[(Logistic model)]\label{exa2logit}
In the logistic regression, we assume
$P(Y=1\mid X)=1-P(Y=0\mid X)=\exp(\alpha_0+\int_0^1X(t)\beta_0(t)\,dt)/
(1+\exp(\alpha_0+\int_0^1X(t)\beta_0(t)\,dt))$. It is easy to see
that $\ell(y;a)=ay-\log(1+\exp(a))$ and
$B(X)=\exp(\alpha_0+\int_0^1 X(t)\beta_0(t)\,dt)/(1+\exp(\alpha
_0+\int_0^1X(t)\beta_0(t)))^2\leq1$. Assumption~\ref{A1}(a) follows from simple algebra. Assumption~\ref{A1}(b)
follows from data transformation and the following $L^2$ bounded
condition: $\int_0^1X^2(t)\,dt\leq c$ a.s. The latter condition implies
that the range\vspace*{1pt} of $\mu_0(X)$ is finite, and thus $B(X)$ is bounded
away from zero. Since $\epsilon=Y-\exp(X^T\theta_0+g_0(Z))/(1+\exp
(X^T\theta_0+g_0(Z)))$, Assumption~\ref{A1}(c) can be
verified by direct calculations.
\end{Example}

%
\begin{Example}[(Exponential family)]\label{exa2exponentialfamily}
Let $(Y,X)$ follow the one-parameter exponential family
\[
Y\mid X\sim\exp\biggl\{Y\biggl(\alpha_0+\int_0^1
X(t)\beta_0(t)\,dt\biggr)+A(Y)-G\biggl(\alpha_0+\int
_0^1 X(t)\beta_0(t)\,dt\biggr)
\biggr\},
\]
where $A(\cdot)$ and $G(\cdot)$ are known, and $\dot{G}=F$ [recall that
$F$ is the link function satisfying (\ref{gflm})]. We assume that $G$
has bounded second- and third-order derivatives,
and $\ddot{G}\ge\delta$ for some constant $\delta>0$; see similar
conditions on page~738\vspace*{1pt} of \cite{MS05}. It is easy to see that $\ell
(y;a)=ya+A(y)-G(a)$, and hence, $\dot{\ell}_a(y;a)=y-\dot{G}(a)$,
$\ddot{\ell}_a(y;a)=-\ddot{G}(a)$ and $\ell'''_a(y;a)=-\dddot
{G}(a)$. Clearly, $\ddot{\ell}_a$ and $\ell'''_a$ are both bounded,
and hence
Assumption~\ref{A1}(a) holds.
Furthermore, $B(X)=\ddot{G}(\alpha_0+\int_0^1 X(t)\beta_0(t)\,dt)$
satisfies Assumption~\ref{A1}(b).
Since $\epsilon=Y-\dot{G}(\alpha_0+\int_0^1 X(t)\beta
_0(t)\,dt)=Y-\mu_0(X)$,
it is easy to see that $E\{\epsilon\mid X\}=E\{Y\mid X\}-\mu_0(X)=0$,
and $E\{\epsilon^2\mid X\}=\operatorname{Var}(Y\mid X)=\ddot{G}(\alpha_0+\int_0^1
X(t)\beta_0(t)\,dt)$
(see \cite{M82}), and therefore, Assumption~\ref{A1}(c) holds.
\end{Example}

\subsection{Reproducing kernel Hilbert space}

We introduce an inner product in $H^m(\mathbb{I})$, that is,
%
\begin{equation}
\label{innerprod1} \langle\beta,\widetilde{\beta}\rangle_1=V(\beta,
\widetilde{\beta})+\lambda J(\beta,\widetilde{\beta}),
\end{equation}
where $V(\beta,\widetilde{\beta})\equiv\int_0^1\int_0^1
C(s,t)\beta(t)\widetilde{\beta}(s)\,ds\,dt$ and $C(s,t)\equiv E\{
B(X)X(t)X(s)\}$ is a weighted covariance function. Denote the
corresponding norm as $\llVert\cdot\rrVert_1$. Define a
linear bounded operator $C(\cdot)$ from $L^2(\mathbb{I})$ to
$L^2(\mathbb{I})$: $(C\beta)(t)=\int_0^1 C(s,t)\beta(s)\,ds$. Below
we assume a regularity condition on $C\beta$, which implies the
positive definiteness of $V$, such that the above inner product (\ref
{innerprod1}) is well defined.

%
\begin{Assumption}\label{A2}
$C(s,t)$ is continuous on $\mathbb{I}\times\mathbb{I}$.
Furthermore, for any $\beta\in L^2(\mathbb{I})$ satisfying $C\beta
=0$, we have $\beta=0$.
\end{Assumption}

Suppose that $C$ is continuous over $\mathbb{I}\times
\mathbb{I}$. By Mercer's theorem, $C$ admits the spectral
decomposition\vspace*{1pt} $C(s,t)=\sum_{\nu=1}^\infty\zeta_\nu\psi_\nu
(s)\psi_\nu(t)$, where $\{\psi_\nu(\cdot), \zeta_\nu\ge0\}_{\nu
\ge1}$ forms an orthonormal basis in $L^2(\mathbb{I})$ under the
usual $L^2$-norm. Therefore, for any $\beta\in L^2(\mathbb I)$, we
have $\beta(\cdot)=\sum_{\nu=1}^\infty b_\nu\psi_\nu(\cdot)$
and $V(\beta,\beta)=\sum_{\nu=1}^\infty\zeta_\nu b_\nu^2$ for a
sequence of square summable $b_\nu$'s. Assumption~\ref{A2} directly
implies that all the eigenvalues of $C$ are positive,
that is, $\zeta_\nu>0$ for all $\nu\ge1$. Therefore, if \mbox{$V(\beta,\beta
)=0$}, that is, $\sum_{\nu=1}^\infty\zeta_\nu b_\nu^2=0$, we
can easily show that $\beta=\sum_{\nu=1}^\infty b_\nu\psi_\nu=0$.
Hence $\langle\cdot,\cdot\rangle_1$ is well defined. Moreover,
together with Proposition 2 of \cite{YC10}, Assumption~\ref{A2}
implies that $H^m(\mathbb{I})$ is indeed a reproducing kernel Hilbert
space (RKHS) under $\langle\cdot,\cdot\rangle_1$. We denote its
reproducing kernel function as $K(s,t)$.

As for the joint parameter space $\mathcal H$, we also need to assume a
proper inner product under which\vspace*{1pt} it is a well-defined Hilbert space.
Define, for any $\theta=(\alpha,\beta)$, $\widetilde{\theta
}=(\widetilde{\alpha},\widetilde{\beta})\in\mathcal{H}$,
%
\begin{eqnarray}\label{innerprod}
\langle\theta, \widetilde{\theta}\rangle&\equiv& E
\biggl\{B(X) \biggl(\alpha+\int_0^1X(t)\beta(t)\,dt \biggr)
\biggl(\widetilde{\alpha}+\int_0^1X(t)
\widetilde{\beta}(t)\,dt \biggr) \biggr\}
\nonumber\\[-8pt]\\[-8pt]\nonumber
&&{} +\lambda J(\widetilde\beta, \beta).\nonumber
\end{eqnarray}
By rewriting $\langle\theta, \widetilde{\theta}\rangle=E\{B(X)\}
\alpha\widetilde\alpha+\langle\beta,\widetilde\beta\rangle_1$,
we note that $\langle\cdot,\cdot\rangle$ is a well-defined inner
product under Assumptions~\ref{A1}(b) and~\ref{A2}. The
corresponding norm is denoted as $\llVert\cdot\rrVert$.
Given the above relation between $\langle\cdot,\cdot\rangle$ and
$\langle\cdot,\cdot\rangle_1$, it is easy to show that $\mathcal H$
inherits the completeness of $H^m(\mathbb I)$. This means $\mathcal
{H}$ is indeed a Hilbert space as described in Proposition~\ref
{Hhilbert} below.

%
\begin{Proposition}\label{Hhilbert}
Under $\langle\cdot,\cdot\rangle$, $\mathcal{H}$ is a Hilbert space.
\end{Proposition}

In the literature, the estimation/prediction rate results in the
(generalized) functional linear models are mostly expressed in terms of
$L^2$-norm; see \cite{MS05,HMW06,HH07,YC10}. We remark that our norm
$\llVert\cdot\rrVert$ is stronger than the $L^2$-norm used
in the above literature under Assumption~\ref{A1}(b).

We next assume a sequence of basis functions in $H^m(\mathbb{I})$
which can simultaneously diagonalize $V$ and $J$. Let $\llVert
\cdot\rrVert_{L^2}$ and $\llVert\cdot\rrVert_{\sup}$
denote the $L^2$ and supremum norms, respectively. Denote $a_n\asymp
b_n$ if and only if there exist positive
constants $c_1,c_2$ such that $c_1\le a_\nu/b_\nu\le c_2$ for all
$\nu$.

%
\begin{Assumption}\label{A3}
There exists a sequence of functions $\{\varphi_\nu\}_{\nu\ge
1}\subset H^m(\mathbb{I})$ such that
$\llVert\varphi_\nu\rrVert_{L^2}\le C_\varphi\nu^{a}$
for each $\nu\ge1$, some constants $a\ge0, C_\varphi>0$ and
%
\begin{equation}
\label{simuldiag} V(\varphi_\nu,\varphi_\mu)=
\delta_{\nu\mu},\qquad J(\varphi_\nu,\varphi_\mu)=
\rho_\nu\delta_{\nu\mu}\qquad\mbox{for any }\nu,\mu\ge1,
\end{equation}
where $\delta_{\nu\mu}$ is Kronecker's notation, and $\rho_\nu$ is
a nondecreasing nonnegative
sequence satisfying $\rho_\nu\asymp\nu^{2k}$ for some constant $k>a+1/2$.
Furthermore, any $\beta\in H^m(\mathbb{I})$ admits the
Fourier expansion $\beta=\sum_{\nu=1}^\infty V(\beta,\varphi_\nu
)\varphi_\nu$ with
convergence in $H^m(\mathbb{I})$ under $\langle\cdot,\cdot\rangle_1$.
\end{Assumption}

We remark that Assumption~\ref{A3} is the price we need to
pay for
making valid statistical inference in addition to those required for
minimax estimation, as in, for example, \cite{YC10,CY12}.

Assumption~\ref{A3} can be directly implied by the pseudo
Sacks--Ylvisaker (SY) conditions, which are slightly different from the
conventional SY conditions proposed in \cite
{SY1966,SY1968,SY1970,RWW95}; see Section~S.2 in
the Supplementary Material \cite{SCFDA} for more details. Proposition
\ref{validA3} below discusses the construction of an eigen-system
satisfying Assumption~\ref{A3} under this condition.

%
\begin{Proposition}[(Eigen-system construction)]\label{validA3}
Suppose the covariance function $C$
satisfies Assumption~\ref{A2} and the pseudo SY conditions of order
$r\ge0$ specified in Section~\textup{S.2}.
Furthermore, the boundary value problem \textup{(S.2)} in Section~\textup{S.3} is regular in the sense of \cite{Birk1908}.
Consider the following integro-differential equations:
%
\begin{equation}
\label{bvp1} \cases{\displaystyle (-1)^m y^{(2m)}_\nu(t)=
\rho_\nu\int_0^1C(s,t)y_\nu(s)\,ds,
\vspace*{5pt}\cr
\displaystyle y_\nu^{(j)}(0)=y_\nu^{(j)}(1)=0, \qquad\mbox{$j=m,\ldots,2m-1$.}}
\end{equation}
Let $(\rho_\nu,y_\nu)$ be the corresponding eigenvalues and
eigenfunctions of problem (\ref{bvp1}),
and let $\varphi_\nu=y_\nu/\sqrt{V(y_\nu,y_\nu)}$. Then $(\rho
_\nu,\varphi_\nu)$ satisfy Assumption~\ref{A3}
with $k=m+r+1$ and $a=r+1$ if \textit{one} of the following additional
assumptions is satisfied:
\begin{longlist}[(ii)]
\item[(i)]$r=0$;
\item[(ii)]$r\geq1$, and for $j=0,1,\ldots,r-1$, $C^{(j,0)}(0,t)=0$ for\vspace*{1pt}
any $0\le t\le1$,
where $C^{(j,0)}(s,t)$ is the $j$th-order partial derivative with
respect to $s$.
\end{longlist}
\end{Proposition}

The proof of Proposition~\ref{validA3} relies on a
nontrivial application of
the general integro-differential equation theory developed in \cite
{T1927,TL1928,T1930}.
In particular, the order of $\rho_\nu$ in problem (\ref{bvp1}) is,
in general, equivalent
to the order of eigenvalues in an ordinary differential problem; see
\cite{T1927}, Theorem 7.
More explicitly, $\rho_\nu\approx(c\pi\nu)^{2k}$ as $\nu
\rightarrow\infty$ for some constant $c>0$; see
\cite{TL1928}, equation (20).

In the Gaussian model with unit variance (see Example~\ref{exa1pss}),
it can be shown with arguments similar to those
in \cite{SY1966,SY1968,SY1970} that
the covariance function $C$ satisfies the pseudo SY conditions of order $r=0$
when $X(t)$ is Brownian motion with $C(s,t)=\min\{s,t\}$.
We also note that the boundary condition in Proposition~\ref{validA3}(ii) was also assumed in \cite{RWW95}
when $X$ is Gaussian of order $r>0$. The integro-differential equations
(\ref{bvp1}) can be
translated into easily computable differential equations. More
specifically, we rewrite $y_\nu(t)$ in
(\ref{bvp1}) as $\ddot{g}_\nu(t)$, and thus obtain that
%
\begin{equation}
\label{bvpexample} \cases{ (-1)^{m+1}g_\nu^{(2m+2)}(t)=
\rho_\nu g_\nu(t),
\vspace*{3pt}\cr
g_\nu^{(j)}(0)=g_\nu^{(j)}(1)=0,\qquad\mbox{$j=m+2,\ldots,2m+1$,}
\vspace*{3pt}\cr
g_\nu(0)=\dot{g}_\nu(1)=0.}
\end{equation}
Note that (\ref{bvp1}) and (\ref{bvpexample}) share the same eigenvalues.
Numerical examinations show that
$\rho_\nu\approx(\pi\nu)^{2(m+1)}$.
The function $g_\nu$'s have closed forms
\[
g_\nu(t)=\operatorname{Re} \Biggl(\sum_{j=1}^{2(m+1)}a_{\nu, j}
\exp\bigl(\rho_\nu^{1/(2(m+1))}z_jt\bigr) \Biggr),\qquad
\nu=1,2,\ldots,
\]
where $\operatorname{Re}(\cdot)$ means the real part of a complex number,
$z_1,\ldots,z_{2(m+1)}$ are the complex (distinct) roots of
$z^{2(m+1)}=(-1)^{m+1}$
and $a_{\nu,1},\ldots,a_{\nu,2(m+1)}$ are complex constant coefficients
determined by the boundary value conditions in (\ref{bvpexample}).
It follows by Proposition~\ref{validA3} that the resultant
$\rho_\nu$ and the corresponding scaled functions $\varphi_\nu
=y_\nu/\sqrt{V(y_\nu,y_\nu)}$
satisfy Assumption~\ref{A3},
where recall that $y_\nu$ is the second-order derivative of $g_\nu$.

In the logistic regression (Example~\ref{exa2logit}) or exponential
family models (Example~\ref{exa2exponentialfamily}),
the approach given in Section~S.5 (Supplementary Material \cite{SCFDA})
can be used to find $(\rho_\nu,\varphi_\nu)$
without verifying the pseudo SY conditions. To do so, we need to
replace the kernel function $C$
by its sample version $C_n(s,t)\equiv n^{-1}\sum_{i=1}^n \widehat
{B}(X_i)X_i(s)X_i(t)$,
where $\widehat{B}(X)$ is the plug-in estimate of $B(X)$.

Recall that $K$ is the reproducing kernel function for $H^m(\mathbb
{I})$ under $\langle\cdot,\cdot\rangle_1$.
For any $t\in\mathbb{I}$, define $K_t(\cdot)=K(t,\cdot)\in
H^m(\mathbb I)$. Under Assumption~\ref{A3}, we may write
$K_t=\sum_{\nu\ge1}a_\nu\varphi_\nu$ for a real sequence $a_\nu$.
Clearly, $\varphi_\nu(t)=\langle K_t,\varphi_\nu\rangle_1=a_\nu
(1+\lambda\rho_\nu)$, for all $\nu\ge1$.
So $K_t=\sum_{\nu\ge1}\frac{\varphi_\nu(t)}{1+\lambda\rho_\nu
}\varphi_\nu$. Define $W_\lambda$ as an operator from $H^m(\mathbb
{I})$ to $H^m(\mathbb{I})$ satisfying $\langle W_\lambda\beta,\widetilde
{\beta}\rangle_1=\lambda J(\beta,\widetilde{\beta})$,
for all $\beta,\widetilde{\beta}\in H^m(\mathbb{I})$. Hence
$W_\lambda$ is linear, nonnegative definite and self-adjoint.
For any $\nu\ge1$, write $W_\lambda\varphi_\nu=\sum_\mu b_\mu
\varphi_\mu$.
Then by Assumption~\ref{A3}, for any $\mu\ge1$, $\lambda\rho_\nu
\delta_{\nu\mu}=
\lambda J(\varphi_\nu,\varphi_\mu)=\langle W_\lambda\varphi_\nu,\varphi
_\mu\rangle_1=b_\mu(1+\lambda\rho_\mu)$.
Therefore, $b_\nu=\lambda\rho_\nu/(1+\lambda\rho_\nu)$ and
$b_\mu=0$ if $\mu\neq\nu$,
which implies $W_\lambda\varphi_\nu=\frac{\lambda\rho_\nu
}{1+\lambda\rho_\nu}\varphi_\nu$.
Thus we have shown the following result.

%
\begin{Proposition}\label{PropRPK} Suppose Assumption~\ref{A3} holds.
For any $t\in\mathbb{I}$,
\[
K_t(\cdot)=\sum_\nu
\frac{\varphi_\nu(t)}{1+\lambda\rho_\nu} \varphi_\nu(\cdot),\vspace*{-2pt}
\]
and for any $\nu\ge1$,\vspace*{-2pt}
\[
(W_\lambda\varphi_\nu) (\cdot)= \frac{\lambda\rho_\nu}{1+\lambda\rho_\nu}
\varphi_\nu(\cdot).
\]
\end{Proposition}

Propositions~\ref{propR} and~\ref{propP} below define two
operators, $R_x$ and $P_\lambda$, that will be used in the Fr\'{e}chet
derivatives of the criterion function $\ell_{n, \lambda}$. We first
define $\tau(x)$ as follows. For any $L^2$ integrable function
$x=x(t)$ and $\beta\in H^m(\mathbb{I})$,
$\mathcal{L}_x(\beta)\equiv\int_0^1 x(t)\beta(t)\,dt$ defines a
linear bounded functional. Then by the Riesz representation theorem,
there exists an element in $H^m(\mathbb{I})$, denoted as $\tau(x)$,
such that $\mathcal{L}_x(\beta)=\langle\tau(x),\beta\rangle_1$
for all $\beta\in H^m(\mathbb{I})$.
If we denote $\tau(x)=\sum_{\nu=1}^\infty x^*_\nu\varphi_\nu$,
then $x^*_\nu(1+\lambda\rho_\nu)=\langle\tau(x),\varphi_\nu
\rangle_1=\int_0^1x(t)\varphi_\nu(t)\,dt\equiv x_\nu$ for\vspace*{1pt} any $\nu
\ge1$.
Thus $\tau(x)=\sum_{\nu=1}^\infty\frac{x_\nu}{1+\lambda\rho_\nu
}\varphi_\nu$.

%
\begin{Proposition}\label{propR}
For any $x\in L^2(\mathbb{I})$,
define $R_x=(E\{B(X)\}^{-1},\tau(x))$.
Then $R_x\in\mathcal{H}$ and $\langle R_x,\theta\rangle=\alpha
+\int_0^1 x(t)\beta(t)\,dt$ for any $\theta=(\alpha,\beta)\in
\mathcal{H}$.
\end{Proposition}

It should be noted that $R_x$ depends on $h$ according to the
definition of $\tau(x)$.

%
\begin{Proposition}\label{propP}
For\vspace*{2pt} any $\theta=(\alpha,\beta)\in\mathcal{H}$, define $P_\lambda
\theta=(0,W_\lambda\beta)$.
Then $P_\lambda\theta\in\mathcal{H}$ and
$\langle P_\lambda\theta,\widetilde{\theta}\rangle=\langle
W_\lambda\beta,\widetilde{\beta}\rangle_1$ for any $\widetilde
{\theta}=(\widetilde{\alpha},\widetilde{\beta})\in\mathcal{H}$.
\end{Proposition}

For notational convenience, denote $\Delta\theta=(\Delta\alpha,\Delta
\beta)$
and $\Delta\theta_j=(\Delta\alpha_j,\Delta\beta_j)$ for
$j=1,2,3$. The Fr\'{e}chet derivative of $\ell_{n,\lambda}(\theta)$
w.r.t. $\theta$ is given by
\[
S_{n,\lambda}(\theta)\Delta\theta\equiv D\ell_{n,\lambda}(\theta)\Delta
\theta=\frac{1}{n}\sum_{i=1}^n\dot{
\ell}_a\bigl(Y_i;\langle R_{X_i},\theta\rangle
\bigr)\langle R_{X_i},\Delta\theta\rangle-\langle P_\lambda
\theta,\Delta\theta\rangle.
\]
The second- and third-order Fr\'{e}chet derivatives of $\ell
_{n,\lambda}(\theta)$ can be shown to be, respectively,
\begin{eqnarray*}
&&DS_{n,\lambda}(\theta)\Delta\theta_1\Delta\theta_2
\\[-2pt]
&&\qquad \equiv D^2\ell_{n,\lambda}(\theta)\Delta\theta_1
\Delta\theta_2
\\[-2pt]
&&\qquad = \frac{1}{n}\sum_{i=1}^n \ddot{
\ell}_a\bigl(Y_i;\langle R_{X_i},\theta\rangle
\bigr)\langle R_{X_i},\Delta\theta_1\rangle\langle
R_{X_i},\Delta\theta_2\rangle-\langle P_\lambda
\Delta\theta_1,\Delta\theta_2\rangle
\end{eqnarray*}
and
\begin{eqnarray*}
&&D^2S_{n,\lambda}(\theta)\Delta\theta_1\Delta
\theta_2\Delta\theta_3
\\[-2pt]
&&\qquad \equiv D^3\ell_{n,\lambda}(\theta)\Delta\theta_1
\Delta\theta_2\Delta\theta_3
\\[-2pt]
&&\qquad =\frac{1}{n}\sum_{i=1}^n
\ell'''_a
\bigl(Y_i;\langle R_{X_i},\theta\rangle\bigr) \langle
R_{X_i},\Delta\theta_1\rangle\langle R_{X_i},
\Delta\theta_2\rangle\langle R_{X_i},\Delta
\theta_3\rangle.
\end{eqnarray*}
Define $S_n(\theta)=\frac{1}{n}\sum_{i=1}^n\dot{\ell
}_a(Y_i;\langle R_{X_i},\theta\rangle)R_{X_i}$,
$S(\theta)=E\{S_n(\theta)\}$ and $S_\lambda(\theta)=E\{S_{n,\lambda
}(\theta)\}$
with expectations taken under the true model.

\section{Bahadur representation for functional data}\label{secbrfd}

In this section, we extend the functional Bahadur representation
originally established in the nonparametric regression framework \cite
{SC13} to the generalized functional linear models. This new technical
tool is fundamentally important in the sense that it provides a unified
treatment for various inference problems.

Denote $h=\lambda^{1/(2k)}$, where $k$ is specified in Assumption~\ref
{A3}. An auxiliary norm is introduced for technical purpose: $\llVert
\theta\rrVert_2=\llvert \alpha\rrvert +\llVert\beta\rrVert
_{L^2}$ for any $\theta=(\alpha,\beta)\in\mathcal{H}$. The
following result gives a useful relationship between the two norms
$\llVert\cdot\rrVert_2$ and $\llVert\cdot\rrVert
$. Recall that $a$ is defined in Assumption~\ref{A3}.

%
\begin{Lemma}\label{crucialrelation}
There exists a constant $\kappa>0$ such that for any $\theta\in
\mathcal{H}$,
$\llVert\theta\rrVert_2\le\kappa h^{-(2a+1)/2}\llVert
\theta\rrVert$.
\end{Lemma}

To obtain an appropriate Bahadur representation for the functional
data, we need the following regularity conditions on $X$. Recall that
$\llVert X\rrVert_{L^2}^2=\int_0^1X^2(t)\,dt$.

%
\begin{Assumption}\label{A4}
There exists a constant $s\in(0,1)$ such that
%
\begin{equation}
\label{A41} E\bigl\{\exp\bigl(s\llVert X\rrVert_{L^2}\bigr)\bigr\}<
\infty.
\end{equation}
Moreover, suppose that there exists a constant $M_0>0$ such that for
any $\beta\in H^m(\mathbb{I})$,
%
\begin{equation}
\label{amomconv} E \biggl\{\biggl\llvert\int_0^1
X(t)\beta(t)\,dt\biggr\rrvert^4 \biggr\}\le M_0 \biggl[E
\biggl\{\biggl\llvert\int_0^1 X(t)\beta(t)\,dt
\biggr\rrvert^2 \biggr\} \biggr]^2.
\end{equation}
\end{Assumption}

It is easy to see that (\ref{A41}) holds for any bounded stochastic
process $X$, that is, $\llVert X\rrVert_{L^2}\le c$ a.s.
for some constant $c>0$.
This applies to Example~\ref{exa2logit} which usually requires $X$
to be almost surely bounded in terms of $L^2$-norm.
Equation~(\ref{A41}) also holds for the Gaussian process
as described in Proposition~\ref{validA4} below.
The result applies to Examples~\ref{exa1pss} and~\ref
{exa2exponentialfamily}
where $X$ can be Gaussian.

%
\begin{Proposition}\label{validA4}
If $X$ is a Gaussian process with square-integrable mean function,
then (\ref{A41}) holds for any $s\in(0,1/4)$.
\end{Proposition}

The fourth moment condition (\ref{amomconv})
is valid for $M_0=3$ when $X$ is a Gaussian process; see \cite{YC10}
for more discussions.
The following result shows that (\ref{amomconv}) actually holds in
more general settings.

%
\begin{Proposition}\label{propvalidfourthmomcond}
Suppose $X(t)=u(t)+\sum_{\nu=1}^\infty\xi_\nu\omega_\nu\psi_\nu(t)$,
where $u(\cdot)\in L^2(\mathbb{I})$ is nonrandom,
$\psi_\nu$ is orthonormal $L^2(\mathbb{I})$-basis,
$\omega_\nu$ is a real square-summable sequence
and $\xi_\nu$ are independent random variables
drawn from some symmetric distribution with finite fourth-order moment.
Then for any $\beta(t)=\sum_{\nu=1}^\infty b_\nu\psi_\nu(t)$
with $b_\nu$ being real square-summable,
(\ref{amomconv}) holds with $M_0=\max\{E\{\xi_\nu^4\}/E\{\xi_\nu
^2\}^2,3\}$.
\end{Proposition}

Lemma~\ref{eplemma} below proves a concentration inequality as a
preliminary step in obtaining the Bahadur representation. Denote
$T=(Y,X)\in\mathcal{T}$ as the data variable. Let $\psi_n(T;\theta
)$ be a function over
$\mathcal{T}\times\mathcal{H}$, which might depend on $n$. Define
\[
H_n(\theta)=\frac{1}{\sqrt{n}}\sum_{i=1}^n
\bigl[\psi_n(T_i;\theta)R_{X_i}-E_T
\bigl\{\psi_n(T;\theta)R_X\bigr\} \bigr],
\]
where $E_T\{\cdot\}$
means the expectation w.r.t. $T$. Define $\mathcal{F}_{p_n}=
\{\theta=(\alpha,\beta)\in\mathcal{H}\dvtx  \llvert \alpha\rrvert \le1,
\llVert
\beta\rrVert_{L^2}\le1, J(\beta,\beta)\le p_n\}$,
where $p_n\ge1$.

%
\begin{Lemma}\label{eplemma}
Suppose Assumptions~\ref{A1} to~\ref{A4} hold.
In addition, $\psi_n(T_i;0)=0$ a.s., there exists a constant $C_\psi
>0$ s.t., and
the following Lipschitz continuity holds:
%
\begin{equation}
\label{lipcond} \bigl\llvert\psi_n(T;\theta)-\psi_n(T;
\widetilde{\theta})\bigr\rrvert\le C_\psi\llVert\theta-\widetilde{
\theta}\rrVert_2\qquad\mbox{for any }\theta,\widetilde{\theta}\in
\mathcal{F}_{p_n}.
\end{equation}
Then as $n\rightarrow\infty$,
\[
\sup_{\theta\in\mathcal{F}_{p_n}} \frac{\llVert H_n(\theta)\rrVert
}{p_n^{1/(4m)}\llVert
\theta\rrVert_2^{\gamma}+n^{-1/2}}=O_P\bigl(
\bigl(h^{-1}\log{\log{n}}\bigr)^{1/2}\bigr),
\]
where $\gamma=1-1/(2m)$.
\end{Lemma}

Our last assumption is concerned with the convergence rate of $\widehat
{\theta}_{n,\lambda}$.
Define $r_n=(nh)^{-1/2}+h^k$. Recall that $k$ is specified in
Assumption~\ref{A3}.

%
\begin{Assumption}\label{A5}
$\llVert\widehat{\theta}_{n,\lambda}-\theta_0\rrVert=O_P(r_n)$.
\end{Assumption}

Proposition~\ref{rateconv} states that Assumptions~\ref{A1} to~\ref
{A4} are actually sufficient to imply the above rate of convergence if
the smoothing parameter is properly chosen. Note that no estimation
consistency is required in Proposition~\ref{rateconv}.

%
\begin{Proposition}\label{rateconv}
Suppose that Assumptions~\ref{A1} to~\ref{A4} hold, and that the
following rate conditions on $h$ (or equivalently, $\lambda$) are satisfied:
%
\begin{eqnarray}\label{ratehrateconv}
h&=&o(1),\nonumber
\\
n^{-1/2}h^{-1}&=&o(1),
\\
n^{-1/2}h^{-(a+1)-\sklvafrac{2k-2a-1}{4m}}(
\log{n})^2(\log\log{n})^{1/2}&=&o(1).\nonumber
\end{eqnarray}
Then Assumption~\ref{A5} is satisfied.
\end{Proposition}

It follows by Proposition~\ref{propvalidfourthmomcond} that
if $X(t)=u(t)+\sum_{\nu=1}^\infty\xi_\nu\omega_\nu\psi_\nu(t)$
with $\xi_\nu$ being independent random variables
following symmetric distribution with bounded support, say, $[-N,N]$,
then (\ref{amomconv}) holds. In this case, $X$ is almost surely
$L^2$ bounded since
$\llVert X\rrVert_{L^2}\le\llVert u\rrVert
_{L^2}+\sqrt{\sum_\nu\xi_\nu^2\omega_\nu^2}
\le\llVert u\rrVert_{L^2}+N\sqrt{\sum_\nu\omega_\nu
^2}$, a.s.
Then Proposition~\ref{rateconv} states that if Assumptions~\ref{A1}
to~\ref{A4} hold
and the smoothing parameter is tuned to satisfy (\ref{ratehrateconv}),
then Assumption~\ref{A5} holds for the above $L^2$ bounded~$X$.

Condition (\ref{ratehrateconv}) is satisfied for a suitable range
of $h$. To illustrate this point, we consider the following simple but
representative case.
Under the setup of Proposition~\ref{validA3}, we have $a=r+1$ and $k=m+r+1$.
Suppose $r=0$, that is, $X$ corresponds to zero-order covariance.\vspace*{1pt} Thus
$a=1$ and $k=m+1$.
Denote $h^\ast\asymp n^{-1/(2k+1)}$, $h^{\ast\ast}\asymp
n^{-2/(4k+1)}$ and $h^{\ast\ast\ast}\asymp n^{-1/(2k)}$.
It can be shown that
when $m>(3+\sqrt{5})/4$, $h^\ast$, $h^{\ast\ast}$ and $h^{\ast\ast
\ast}$
all satisfy the conditions of (\ref{ratehrateconv}).
It should be mentioned that $h^\ast$ yields the optimal estimation\vspace*{1pt}
rate $n^{-k/(2k+1)}$ \cite{YC10},
$h^{\ast\ast}$ yields the optimal testing rate $n^{-2k/(4k+1)}$ as
will be shown in later sections,
and $h^{\ast\ast\ast}$ yields the optimal prediction rate \cite{CH06}.

Now we are ready to present the Bahadur representation based on the
functional data.

\begin{Theorem}[(Bahadur representation for functional data)] \label{UBR}
Suppose that Assumptions
\ref{A1}--\ref{A5} hold, and as $n\rightarrow\infty$,
$h=o(1)$ and $\log(h^{-1})=O(\log{n})$.
Furthermore, (\ref{amomconv}) holds.
Then, as $n\rightarrow\infty$,
$\llVert\widehat{\theta}_{n,\lambda}-\theta_0-S_{n,\lambda
}(\theta_0)\rrVert=O_P(a_n)$,
where
\[
a_n=n^{-1/2}h^{-\vafrac{4ma+6m-1}{4m}}r_n(
\log{n})^2(\log\log{n})^{1/2}+C_\ell
h^{-1/2}r_n^2,
\]
and
\[
C_\ell\equiv\sup_{x\in L^2(\mathbb{I})}E\Bigl\{\sup
_{a\in\bbR}\bigl\llvert\ell'''_a(Y;a)
\bigr\rrvert \big| X=x\Bigr\}.
\]
\end{Theorem}

We next give an example rate of $a_n$ in Theorem~\ref{UBR} when $\ell
$ is quadratic. In this case, we have $C_\ell=0$. Suppose $a=1$ and
$k=m+1$; see the discussions below Proposition~\ref{rateconv}.
Direct examinations show that $a_n$ is of the order $o(n^{-1/2})$ when
$m>1+\sqrt{3}/2\approx1.866$
and $h=h^\ast$, $h^{\ast\ast}$ and $h^{\ast\ast\ast}$.

An immediate consequence of Bahadur representation is the following
point-wise limit distribution of the slope function estimate. This
local result is new and of independent interest, for example,
point-wise CI.

%
\begin{Corollary}\label{localasymp}
Suppose that the conditions of Theorem~\ref{UBR} are satisfied, $\sup
_{\nu\ge1}\llVert\varphi_\nu\rrVert_{\sup}
\le C_{\varphi}\nu^a$ for $\nu\ge1$ and that
$E\{\exp(s\llvert \epsilon\rrvert )\}<\infty$, for some constant $s>0$.
Furthermore, as $n\rightarrow\infty$, $nh^{2a+1}(\log
(1/h))^{-4}\rightarrow\infty$,
$n^{1/2}a_n=o(1)$ and $\sum_{\nu=1}^\infty
\frac{\llvert \varphi_\nu(z)\rrvert ^2}{(1+\lambda\rho_\nu)^2}\asymp
h^{-(2a+1)}$.
Then we have for any $z\in\mathbb{I}$,
\[
\frac{\sqrt{n}(\widehat{\beta}_{n,\lambda}(z)-\beta
_0(z)+(W_\lambda\beta_0)(z))}{
\sqrt{\sum_{\nu=1}^\infty
\sklfrac{\llvert \varphi_\nu(z)\rrvert ^2}{(1+\lambda\rho_\nu)^2}}}
\stackrel{d} {\longrightarrow}N(0,1).
\]
In addition, if $\sqrt{n}(W_\lambda\beta
_0)(z) /\sqrt{\sum_{\nu
=1}^\infty
\frac{\llvert \varphi_\nu(z)\rrvert ^2}{(1+\lambda\rho_\nu)^2}}=o(1)$,
then
\[
\frac{\sqrt{n}(\widehat{\beta}_{n,\lambda}(z)-\beta_0(z))}{
\sqrt{\sum_{\nu=1}^\infty
\sfrac{\llvert \varphi_\nu(z)\rrvert ^2}{(1+\lambda\rho_\nu)^2}}}
\stackrel{d} {\longrightarrow}N(0,1).
\]
\end{Corollary}

Corollary~\ref{localasymp} applies\vspace*{1pt} to any point $z\in\mathbb{I}$ satisfying
$\sum_{\nu=1}^\infty
\frac{\llvert \varphi_\nu(z)\rrvert ^2}{(1+\lambda\rho_\nu)^2}\asymp
h^{-(2a+1)}$.
Validity of this condition is discussed in Section~S.14
of the Supplementary Material~\cite{SCFDA}.

The condition $\sqrt{n}(W_\lambda\beta
_0)(z)/\sqrt{\sum_{\nu
=1}^\infty
\frac{\llvert \varphi_\nu(z)\rrvert ^2}{(1+\lambda\rho_\nu)^2}}=o(1)$
holds if
$nh^{4k}=o(1)$ and the true slope function $\beta_0=\sum_\nu b_\nu
\varphi_\nu$ satisfies the condition (\textbf{U}):\break $\sum_\nu b_\nu
^2\rho_\nu^2<\infty$. To see this,
observe that
\begin{eqnarray*}
\bigl\llvert(W_\lambda\beta_0) (z)\bigr\rrvert
&=&\biggl\llvert\sum_\nu b_\nu
\frac{\lambda\rho_\nu}{1+\lambda\rho_\nu
}\varphi_\nu(z)\biggr\rrvert
\\
&\le& C_\varphi\lambda\sum_{\nu}\llvert
b_\nu\rrvert\frac{\rho_\nu\nu
^a}{1+\lambda\rho_\nu}
\\
&\le& C_\varphi\lambda\sqrt{\sum_\nu
b_\nu^2\rho_\nu^2} \sqrt{\sum
_\nu\frac{\nu^{2a}}{(1+\lambda\rho_\nu)^2}},
\end{eqnarray*}
where the last term is of the order $O(\lambda h^{-(2a+1)/2})$. Hence,
it leads to
\[
\sqrt{n}(W_\lambda\beta_0) (z)\bigg/\sqrt{\sum
_{\nu=1}^\infty\frac{\llvert \varphi_\nu(z)\rrvert ^2}{(1+\lambda\rho
_\nu)^2}}\asymp\sqrt
{nh^{2a+1}}(W_\lambda\beta_0) (z)=o(1).
\]

%
\begin{Remark}[(Convergence rate)]
Corollary~\ref{localasymp} derives the convergence
rate of the local estimate $\widehat{\beta}_{n,\lambda}(z)$ as
$\sqrt{nh^{2a+1}}$.
The factor $a$ (defined in Assumption~\ref{A3}) generically reflects
the impact of the covariance operator on the convergence rate. For
example, Proposition~\ref{validA3} shows that $a=r+1$ with $r$ being
the order of the covariance function under the pseudo SY condition. The
above observation coincides with the arguments in \cite{HH07}, that
the covariance effect in general influences the (global) rate
expressions. When the eigenfunctions are uniformly bounded, that is,
$a=0$, the above rate becomes $\sqrt{nh}$, which is exactly the rate
derived in the general nonparametric regression setup; see Theorem 3.5
in \cite{SC13}.
\end{Remark}

%
\begin{Remark}
(Undersmoothing) Assumption~\ref{A3} implies that $\beta_0\in
H^m(\mathbb{I})$ has the property that $\sum_\nu b_\nu^2\rho_\nu
<\infty$. However, the condition (\textbf{U}) imposes a faster decay
rate on the generalized Fourier coefficients $b_\nu$, which in turn
requires more smoothness of $\beta_0$. Since we still employ the $m$th
order penalty in~(\ref{estimation}), condition~(\textbf{U}) can be
treated as a type of undersmoothing condition. More generally, similar
conditions will be implicitly imposed in the inference procedures to be
presented later.
\end{Remark}

\section{Confidence/prediction interval}\label{secasympCI}
In this section, we consider two inter-connected inference procedures:
(i) confidence interval for the conditional mean and (ii) prediction
interval for a future response.

\subsection{Confidence interval for conditional mean}\label{secCI}
For any (nonrandom) $x_0\in L^2(\mathbb{I})$,
we construct a confidence interval $\mu_0(x_0)=E\{Y\mid X=x_0\}$ by
centering around the plug-in estimate $\widehat{Y}_0\equiv F(\widehat
{\alpha}_{n,\lambda}+\int_0^1 x_0(t)\widehat{\beta}_{n,\lambda}(t)\,dt)$.
Define $\mu_0'(x_0)=\dot{F}(\alpha_0+\int_0^1x_0(t)\beta_0(t)\,dt)$,
$\sigma_n^2=E\{B(X)\}^{-1}+\sum_{\nu=1}^\infty\frac{\llvert x_\nu
^0\rrvert ^2}{(1+\lambda\rho_\nu)^2}$,
where $x_\nu^0=\int_0^1 x_0(t)\varphi_\nu(t)\,dt$.

%
\begin{Theorem}[(Confidence interval construction)]\label{CImean}
Let Assumptions~\ref{A1} through
\ref{A5} be satisfied for the true parameter
$\theta_0=(\alpha_0,\beta_0)$,
and $\mu_0'(x_0)\neq0$. Furthermore,\vspace*{1pt} assume (\ref{amomconv}) and
$E\{\exp(s\llvert \epsilon\rrvert )\}<\infty$ for some $s>0$.
If $h=o(1)$, $\log(h^{-1})=O(\log{n})$, $nh^{2a+1}(\log{n})^{-4}
\rightarrow\infty$, $na_n^2=o(1)$ and $\llVert R_{x_0}\rrVert\asymp
\sigma_n$,
then as $n\rightarrow\infty$,
\begin{eqnarray*}
&&\frac{\sqrt{n}}{\sigma_n} \biggl(\widehat{\alpha}_{n,\lambda
}+\int
_0^1 x_0(t)\widehat{
\beta}_{n,\lambda}(t)\,dt -\alpha_0-\int_0^1x_0(t)
\beta_0(t)\,dt
\\
&&\hspace*{137pt}{} -\int_0^1
x_0(t) (W_\lambda\beta_0) (t)\,dt \biggr)
\\[-1pt]
&&\qquad\stackrel{d} {\longrightarrow}N(0,1).
\end{eqnarray*}

Furthermore, if $\beta_0=\sum_\nu b_\nu\varphi_\nu$
with $\sum_\nu b_\nu^2\rho_\nu^2<\infty$ and $nh^{4k}=o(1)$, then
$\frac{\sqrt{n}}{\sigma_n}\int_0^1 x_0(t)(W_\lambda\beta
_0)(t)\,dt=o(1)$ so that
we have
%
\begin{equation}
\label{CIlimitdistr} \frac{\sqrt{n}}{\sigma_n\mu_0'(x_0)} \bigl
(\widehat Y_0-\mu
_0(x_0) \bigr)\stackrel{d} {\longrightarrow}N(0,1).
\end{equation}
Hence the $100(1-\widetilde{\alpha})\%$ confidence interval for $\mu
_0(x_0)$ is
%
\begin{equation}
\label{CIconditionalmean} \bigl[\widehat Y_0 \pm n^{-1/2}z_{\widetilde
{\alpha}/2}
\sigma_n \widehat\mu_0'(x_0)
\bigr],
\end{equation}
where $z_{\widetilde{\alpha}/2}$ is the $(1-\widetilde{\alpha
}/2)$-quantile of $N(0,1)$ and $\widehat\mu_0'(x_0)\equiv\dot
{F}(\widehat\alpha_{n, \lambda}+\break \int_0^1x_0(t)\widehat\beta
_{n,\lambda}(t)\,dt)$.
\end{Theorem}

In\vspace*{1pt} the Gaussian model (Example~\ref{exa1pss}) with $B(X)\equiv1$,
if $X$ is Brownian motion with $C(s,t)=\min\{s,t\}$,
then $\sigma_n^2$ has an explicit form with $\rho_\nu\approx(2\pi
\nu)^{2(m+1)}$
and $\varphi_\nu$ solved by (\ref{bvpexample}). As for Examples
\ref{exa2logit} and~\ref{exa2exponentialfamily}, one can obtain
$\sigma_n^2$ by following the approach outlined in Section~S.5
(Supplementary Material \cite{SCFDA}).

A direct byproduct of Theorem~\ref{CImean} is the prediction rate
$\sigma_n/\sqrt{n}$. Proposition~\ref{proppredrate} further
characterizes this rate in various situations. Suppose that $\llvert
x_\nu
^0\rrvert \asymp\nu^{a-d}$ for some constant $d$.
A larger $d$ usually yields a smoother function $x_0$.\vadjust{\goodbreak}

%
\begin{Proposition}\label{proppredrate}
The prediction rate in Theorem~\ref{CImean} satisfies
\[
\sigma_n/\sqrt{n}= \cases{ n^{-1/2}, &\quad if $d-a>1/2$,
\vspace*{3pt}\cr
n^{-1/2}\bigl(\log(1/h)\bigr)^{1/2}, &\quad if $d-a=1/2$,
\vspace*{3pt}\cr
n^{-1/2}h^{d-a-1/2}, &\quad if $d-a<1/2$.}
\]
In\vspace*{1pt} particular, if $d-a<1/2$ and $h=h^{\ast\ast\ast}\asymp n^{-\afrac
{1}{2(m+a)}}$, then
$\sigma_n/\sqrt{n}=n^{-\vafrac{d+m-1/2}{2(m+a)}}$. Furthermore, if $k=m+a$
as in the setting of Proposition~\ref{validA3}, then $\sigma_n/\sqrt
{n}$ is minimax optimal when $h=h^{\ast\ast\ast}$.
\end{Proposition}

Proposition~\ref{proppredrate} states that when $d-a>1/2$, that is,
the process $x_0$ is sufficiently smooth, then the prediction can be
conducted in
terms of root-$n$ rate regardless of the choice of $h$. This result
coincides with \cite{CH06}
in the special FPCA setting. Moreover, when $d-a<1/2$ and $h=h^{\ast
\ast\ast}$, the rate becomes optimal. Again, this is consistent with
\cite{CH06} in the setting
that the true slope function belongs to a Sobolev rectangle.
Interestingly, it can be checked that $h=h^{\ast\ast\ast}$
satisfies the rate conditions in Theorem~\ref{CImean}
if $a=1$, $k=m+1$ and $m>1+\sqrt{3}/2$; see the discussions below
Theorem~\ref{UBR}.

\subsection{Prediction interval for future response}

Following Theorem~\ref{CImean}, we can establish the prediction
interval for the future response $Y_0$ conditional on \mbox{$X=x_0$}. Write
$Y_0-\widehat Y_0=\xi_n+\epsilon_0$, where $\xi_n=\mu
_0(x_0)-\widehat Y_0$ and $\epsilon_0=\dot\ell_a(Y_0;\alpha_0+\int
_0^1x_0(t)\beta_0(t)\,dt)$. Since $\epsilon_0$ is independent of $\xi
_n$ depending on all the past data $\{Y_i,X_i\}_{i=1}^n$, we can easily
incorporate the additional randomness from $\epsilon_0$ into the
construction of the prediction interval. This leads to a nonvanishing
interval length as sample size increases. This is crucially different
from that of confidence interval.

Let $F_{\xi_n}$ and $F_{\epsilon_0}$ be the distribution functions of
$\xi_n$ and $\epsilon_0$, respectively.
Denote the distribution function of $\xi_n+\epsilon_0$ as $G\equiv
F_{\xi_n}*F_{\epsilon_0}$, and $(l_{\widetilde{\alpha}},
u_{\widetilde{\alpha}})$ as its $(\widetilde{\alpha}/2)$th and
$(1-\widetilde{\alpha}/2)$th quantiles, respectively. Then the
$100(1-\widetilde{\alpha})\%$ prediction interval for $Y_0$ is given as
\[
[\widehat Y_0+l_{\widetilde{\alpha}}, \widehat Y_0+u_{\widetilde{\alpha}}
].
\]
Theorem~\ref{CImean} directly implies that $\xi_n\stackrel{a}{\sim
}N(0, (n^{-1/2}\sigma_n\widehat\mu_0'(x_0))^2)$,
where $\stackrel{a}{\sim}$ means \textit{approximately distributed}. If
we further assume that $\epsilon_0\sim N(0, B^{-1}(x_0))$ [see
Assumption~\ref{A1}(c)], that is, $B(x_0)$ is the reciprocal
error variance for the $L^2$ loss, the above general formula reduces to
%
\begin{equation}
\label{PIGaussian} \Bigl[\widehat Y_0\pm z_{\widetilde{\alpha}/2}\sqrt
{B(x_0)+\bigl(n^{-1/2}\sigma_n\widehat
\mu_0'(x_0)\bigr)^2} \Bigr].
\end{equation}
The unknown quantities in (\ref{CIconditionalmean}) and (\ref
{PIGaussian}) can be estimated by plug-in approach.

\section{Hypothesis testing}\label{sechypothesistesting}
We consider two types of testing for the generalized functional linear
models: (i) testing the \textit{functional contrast} defined as $\int_0^1
w(t)\beta(t)\,dt$ for some given weight function $w(\cdot)$, for
example, $w=X$ and (ii) testing the intercept value and the global
behavior of the slope function, for example, $\alpha=0$ and $\beta$
is a linear function.

\subsection{Testing functional contrast}

In practice, it is often of interest to test the \emph{functional
contrast}. For example, we may test single frequency or frequency
contrast of the slope function; see Examples~\ref{egCT1} and~\ref
{egCT2}. In general, we test $H_0^{CT}\dvtx \int_0^1 w(t)\beta(t)\,dt=c$
for some known $w(\cdot)$ and $c$.

Consider the following test statistic:
%
\begin{equation}
\label{testcontr} CT_{n,\lambda}=\frac{\sqrt{n}(\int_0^1w(t)\widehat
{\beta
}_{n,\lambda}(t)\,dt-c)}{
\sqrt{\sum_{\nu=1}^\infty\sklfrac{w_\nu^2}{(1+\lambda\rho_\nu)^2}}},
\end{equation}
where $w_\nu=\int_0^1 w(t)\varphi_\nu(t)\,dt$. Recall that $(\varphi
_\nu,\rho_\nu)$ is the eigensystem satisfying Assumption~\ref{A3}.
Let $w\in L^2(\mathbb{I})$ and $\tau(w)\in H^m(\mathbb{I})$ be such
that $\langle\tau(w),\beta\rangle_1=\int_0^1w(t)\beta(t)\,dt$,
for any $\beta\in H^m(\mathbb{I})$. We can verify that $\tau(w)=\sum
_{\nu=1}^\infty\frac{w_\nu}{1+\lambda\rho_\nu}\varphi_\nu$.
Then, under $H_0^{CT}$, $CT_{n,\lambda}$ can be rewritten as
%
\begin{equation}
\frac{\sqrt{n}\langle\tau(w), (\widehat\beta_{n,\lambda}-\beta
)\rangle_1}{\llVert\tau(w)\rrVert_1}.\label{cttest}
\end{equation}
It follows from Theorem~\ref{CImean} that
(\ref{cttest}) converges weakly to a standard normal distribution.
This is summarized in the following theorem.
Define $M_a=\sum_{\nu=1}^\infty\frac{w_\nu^2}{(1+\lambda\rho_\nu
)^a}$ for $a=1,2$.

%
\begin{Theorem}[(Functional contrast testing)]\label{limitCT}
Suppose that Assumptions~\ref{A1}
through~\ref{A5} hold. Furthermore,
let $\beta_0=\sum_\nu b_\nu\varphi_\nu$ with
$\sum_\nu b_\nu^2\rho_\nu^2<\infty$, and
assume~(\ref{amomconv}), $E\{\exp(s\llvert \epsilon\rrvert )\}<\infty
$ for some
$s>0$, and as $n\rightarrow\infty$,
$h=o(1)$, $\log(h^{-1})=O(\log{n})$,
$nh^{2a+1}(\log(1/h))^{-4}\rightarrow\infty$,
$M_1\asymp M_2$, $na_n^2=o(1)$, $nh^{4k}=o(1)$.
Then, under $H_0^{CT}$, we have
$CT_{n,\lambda}\stackrel{d}{\longrightarrow}N(0,1)$ as $n\rightarrow
\infty$.
\end{Theorem}

%
\begin{Example}[(Testing single frequency)]\label{egCT1}
Suppose that the slope function has an expansion $\beta=\sum_{\nu
=1}^\infty b_\nu\varphi_\nu$,
and we want to test whether $b_{\nu^\ast}=0$ for some \mbox{$\nu^\ast\ge
1$}. In other words,
we are interested in knowing whether the \mbox{$\nu^\ast$-}level frequency
of $\beta$ vanishes.
Let $w(t)=(C\varphi_{\nu^\ast})(t)$. Then it is easy to see that
$\int_0^1w(t)\*\beta(t)\,dt=b_{\nu^\ast}$. That is, the problem reduces
to testing
$H_0^{CT}\dvtx \int_0^1w(t)\beta(t)\,dt=0$.
It can be shown directly that $M_a=(1+\lambda\rho_{\nu^\ast
})^{-a}\asymp1$ for $a=1,2$.
If $r=0$ (see Proposition~\ref{validA4} for validity),
then it can be shown that when $m>(3+\sqrt{5})/4\approx1.309$, the
rate conditions in Theorem~\ref{limitCT}
are satisfied for $h=h^\ast$. This means that $H_0^{CT}$ is rejected
at level $0.05$
if $\llvert CT_{n,\lambda}\rrvert >1.96$.
\end{Example}

%
\begin{Example}[(Testing frequency contrast)]\label{egCT2}
Following Example~\ref{egCT1}, we now test whether $\sum_{\nu
=1}^\infty\mathfrak{c}_\nu b_\nu=0$,
for some real sequence $\mathfrak{c}_\nu$ satisfying $0<\inf_{\nu
\ge1}\llvert \mathfrak{c}_\nu\rrvert \le
\sup_{\nu\ge1}\llvert \mathfrak{c}_\nu\rrvert <\infty$.
Suppose that the covariance function $C(\cdot,\cdot)$ satisfies the
conditions in Proposition
\ref{validA3} with order $r>0$.
It follows from\vspace*{1pt} Proposition~\ref{validA3} and its proof that the
eigenfunction $\varphi_\nu$
can be managed so that $\llVert\varphi^{(2m)}_\nu\rrVert
_{\sup}\le C_\varphi\nu^{2m+r+1}$
and $\llVert C\varphi_\nu\rrVert_{\sup}\le\rho_\nu
^{-1}\llVert\varphi^{(2m)}_\nu\rrVert_{\sup}\asymp\nu
^{-(r+1)}$, for all $\nu\ge1$.
So the function $w(t)=\sum_{\nu=1}^\infty\mathfrak{c}_\nu(C\varphi
_\nu)(t)$ is well defined since
the series is absolutely convergent on $[0,1]$. It is easy to see that
$\int_0^1w(t)\beta(t)\,dt=\sum_\nu\mathfrak{c}_\nu b_\nu$
and $\mathfrak{c}_\nu=\int_0^1w(t)\varphi_\nu(t)\,dt=w_\nu$.
So the problem reduces to testing $H_0^{CT}\dvtx \int_0^1w(t)\beta(t)\,dt=0$.
For $a=1,2$
\[
M_a=\sum_{\nu=1}^\infty
\frac{w_\nu^2}{(1+\lambda\rho_\nu
)^a}=\sum_{\nu=1}^\infty
\frac{\mathfrak{c}_\nu^2}{(1+\lambda\rho
_\nu)^a} \asymp\sum_{\nu=1}^\infty
\frac{1}{(1+\lambda\rho_\nu)^a}\asymp h^{-1}.
\]
It can be shown that when $m>(3+\sqrt{5})$, the rate conditions in
Theorem~\ref{limitCT}
are satisfied for $h=h^\ast$. We reject $H_0^{CT}$ at level 0.05 if
$\llvert CT_{n,\lambda}\rrvert >1.96$.
\end{Example}

\subsection{Likelihood ratio testing}\label{secPLRT}
Consider the following simple hypothesis:
%
\begin{equation}
\label{simplehypothesis} H_0\dvtx  \theta=\theta_0\quad\mbox{versus}
\quad H_1\dvtx  \theta\in\mathcal H-\{\theta_0\},
\end{equation}
where $\theta_0\in\mathcal{H}$. The penalized likelihood ratio test
statistic is defined as
%
\begin{equation}
\mathrm{PLRT}=\ell_{n,\lambda}(\theta_0)-\ell_{n,\lambda}(\widehat{
\theta}_{n,\lambda}).
\end{equation}
Recall that $\widehat{\theta}_{n,\lambda}$ is the maximizer of $\ell
_{n,\lambda}(\theta)$ over $\mathcal{H}$. The proposed likelihood
ratio testing also applies to the composite hypothesis; that is,
$\theta$ belongs to a certain class. See Remark~\ref{remcomp} for
more details.

Theorem~\ref{FDAPLRT} below derives the null limiting distribution of
$\mathrm{PLRT}_{n,\lambda}$.

\begin{Theorem}[(Likelihood ratio testing)]\label{FDAPLRT}
Suppose that $H_0$ holds, and Assumptions
\ref{A1} through~\ref{A5} are satisfied for
the hypothesized value $\theta_0$. Let $h$
satisfy the following rate conditions: as $n\rightarrow\infty$,
$nh^{2k+1}=O(1)$, $nh\rightarrow\infty$, $n^{1/2}a_n=o(1)$,
$nr_n^3=o(1)$, $n^{1/2}h^{-(a+\sfrac{1}{2}+\vafrac
{2k-2a-1}{4m})}r_n^2(\log{n})^2\times\break (\log\log{n})^{1/2}=o(1)$
and $n^{1/2}h^{-(2a+1+\vafrac{2k-2a-1}{4m})}\times
r_n^3(\log{n})^3\times\break (\log\log{n})^{1/2}=o(1)$.
Furthermore, there exists a constant $M_4>0$ s.t.\break 
$E\{\epsilon^4\mid X\}\le M_4$, a.s.
Then as $n\rightarrow\infty$,
%
\begin{equation}
\label{PLRTlim} -(2u_n)^{-1/2} \bigl(2n\sigma^2
\cdot \mathrm{PLRT}+u_n+n\sigma^2\llVert W_\lambda
\beta_0\rrVert_1^2 \bigr)\stackrel{d} {
\longrightarrow}N(0,1),
\end{equation}
where $u_n=h^{-1}\sigma_1^4/\sigma_2^2$, $\sigma^2=\sigma
_1^2/\sigma_2^2$ and $\sigma_l^2=h\sum_\nu(1+\lambda\rho_\nu
)^{-l}$ for $l=1,2$.
\end{Theorem}

By carefully examining the proof of Theorem~\ref{FDAPLRT},
it can be shown that $n\llVert W_\lambda\beta_0\rrVert
_1^2=o(n\lambda)=o(u_n)$. Therefore, $-2n\sigma^2\cdot \mathrm{PLRT}$ is
asymptotically\break $N(u_n, 2u_n)$ which is nearly $\chi_{u_n}^2$ as
$n\rightarrow\infty$. Hence we claim the null limit distribution as
being approximately $\chi_{u_n}^2$, denoted as
%
\begin{equation}
\label{wilks} -2n\sigma^2\cdot \mathrm{PLRT}\stackrel{a} {\sim}
\chi_{u_n}^2,
\end{equation}
where $\stackrel{a}{\sim}$ means \textit{approximately distributed}; see
\cite{FZZ01}.
If $C$ satisfies the conditions of Proposition~\ref{validA3} with
order $r\ge0$,
then $\rho_\nu\approx(c\nu)^{2k}$
(see the comments below Proposition~\ref{validA3}),
where $k=m+r+1$ and $c>0$ is constant.
It is easy to see that $\sigma_l^2\approx c^{-1}\int_0^\infty
(1+x^{2k})^{-l}\,dx$ for $l=1,2$.
In Example~\ref{exa1pss}, since the covariance function is free of
the model parameters $\alpha_0,\beta_0$, we can see that $c$ is also
free of the model parameters. In particular, when $B(X)\equiv1$, $m=2$
and $X(t)$ is Brownian motion with $C(s,t)=\min\{s,t\}$ and $r=0$, we
have $k=3$ and $c=\pi$ [by solving~(\ref{bvpexample})]. This yields
$\sigma_l^2\approx0.2876697$, $0.2662496$ for $l=1,2$, respectively.
In the end, we obtain $\sigma^2\approx1.080451$ and $u_n\approx
0.3108129/h$ in (\ref{wilks}). As seen above, the null limiting
distribution has the nice property that it is free of the unknown model
parameters, that is, so-called Wilks phenomenon \cite{W38,FZZ01}.
Hence we unveil a new version of Wilks phenomenon that applies to the
functional data. This Wilks type of result enables us to simulate the
null limit distribution directly without resorting to bootstrap or
other resampling methods.

The quantities $\sigma_1^2,\sigma_2^2$
in Theorem~\ref{FDAPLRT} depend on the population eigenpairs.
However, it is possible to replace these quantities by suitable estimators
so that the results become more applicable. In
Section~S.18 (Supplementary Material \cite{SCFDA}), we discuss the validity of this ``plug-in''
approach for Theorems~\ref{CImean},~\ref{limitCT} and~\ref{FDAPLRT}.

%
\begin{Remark}[(Composite hypothesis)]\label{remcomp}
By examining the proof of Theorem~\ref{FDAPLRT}, we find that the
null limiting distribution derived therein
remains the same even when the hypothesized value $\theta_0$ is
unknown. An important consequence is that the
proposed likelihood ratio approach can also be used to test a composite
hypothesis such as $H_0\dvtx \alpha=\alpha_0$ and $\beta\in
\mathcal{P}_j$, where $\mathcal{P}_j$ represents the class of the
$j$th-order polynomials. Under $H_0$, $\beta
$ is of the form $\beta(t)=\sum_{l=0}^j b_l t^l$
for some unknown vector $\mathbf{b}=(b_0,b_1,\ldots,b_j)^T$. In this
case, the slope function and intercept
can be estimated through the following ``parametric'' optimization:
%
\begin{eqnarray}\label{H0mle}
\bigl(\widehat{\alpha}^0,\widehat{
\mathbf{b}}^0\bigr) &=&\arg\max_{\alpha,b_0,\ldots,b_j\in\bbR}
n^{-1}\sum_{i=1}^n\ell
\Biggl(Y_i;\alpha+\sum_{l=0}^j
b_l\int_0^1 X_i(t)t^l\,dt
\Biggr)
\nonumber\\[-8pt]\\[-8pt]\nonumber
&&{}  -(\lambda/2)\mathbf{b}^T D\mathbf{b},
\end{eqnarray}
where $D=[D_{l_1L^2}]_{l_1,l_2=0,\ldots,j}$ is a $(j+1)\times(j+1)$
matrix with $D_{l_1l_2}
=J(t^{l_1},t^{l_2})$.
The corresponding slope function estimate is $\widehat{\beta
}^0(t)=\sum_{l=0}^j\widehat{b}^0_l t^l$.
The test\vspace*{1pt} statistic for this composite hypothesis is defined as
$\mathrm{PLRT}=\ell_{n,\lambda}(\widehat{\alpha}^0,
\widehat{\beta}^0)
-\ell_{n,\lambda}(\widehat{\alpha}_{n,\lambda},\widehat{\beta
}_{n,\lambda})$. Let $\theta_0=(\alpha_0,\beta_0)
$ be the unknown true model parameter under $H_0$, where $\beta_0$ can
be represented as $\sum_{l=0}^j b_l^0
t^l$. Hence, we can further decompose the above $\mathrm{PLRT}$ as
$\mathrm{PLRT}_1-\mathrm{PLRT}_2$, where $\mathrm{PLRT}_1=\ell_{n,\lambda}
(\theta_0)-\ell_{n,\lambda}(\widehat{\alpha}_{n,\lambda},\widehat
{\beta}_{n,\lambda})$,
$\mathrm{PLRT}_2=\ell_{n,\lambda}(\theta_0)-\ell_{n,\lambda}(\widehat
{\alpha}^0,\widehat{\beta}^0)$.
Note that $\mathrm{PLRT}_1$ is the test statistic for the simple hypothesis
$\theta=\theta_0$ versus $\theta\neq
\theta_0$, and $\mathrm{PLRT}_2$ for the parametric hypothesis $(\alpha,\mathbf{b})=(\alpha_0,
\mathbf{b}^0)$ versus $
(\alpha,\mathbf{b})\neq(\alpha_0,\mathbf{b}^0)$, where $\mathbf{b}^0=(b_0^0,\ldots,b_j^0)^T$.
Conventional\vspace*{1pt} theory on parametric likelihood ratio testing leads to
$-2n\cdot \mathrm{PLRT}_2=O_P(1)$. On the other hand, Theorem~\ref{FDAPLRT}
shows that $-2n\sigma^2\cdot
\mathrm{PLRT}_1\stackrel{a}{\sim}\chi_{u_n}^2$. Therefore, we conclude that
the null limit distribution for testing the
composite hypothesis also follows $\chi^2_{u_n}$.
\end{Remark}

%

In the end of this section, we show that the proposed PLRT is optimal
in the minimax sense \cite{I93} when $h=h^{**}$. To derive the
minimax rate of testing (also called as minimum separation rate), we consider
a local alternative written as $H_{1n}\dvtx  \theta=\theta_{n0}$, where
the alternative value is assumed to deviate
from the null value by an amount of $\theta_n$, that is, $\theta
_{n0}=\theta_0+\theta_n$. For simplicity, we
assume $\theta_0=0$, and thus $\theta_{n0}=\theta_n$. Define the
alternative value set $\theta_n\in\Theta_b
\equiv\{(\alpha,\beta)\in\mathcal{H}\dvtx  \llvert \alpha\rrvert \le b,
\llVert
\beta\rrVert_{L^2}\le b, J(\beta,\beta)\le b\}$ for some
fixed constant $b>0$.

%
\begin{Theorem}\label{mrt}
Let Assumptions~\ref{A1}--\ref{A5} be satisfied uniformly under
$\theta=\theta_{n0}\in\Theta_b$.
Let $h$ satisfy $nh^{3/2}\rightarrow\infty$, as $n\rightarrow\infty$,
and also the rate conditions specified in Theorem~\ref{FDAPLRT}.
Furthermore, $\inf_{y\in\mathcal{Y},a\in\bbR}(-\ddot{\ell}_a(y;a))>0$,
and there is a constant $M_4>0$ s.t. for $\theta_n=(\alpha_n,\beta
_n)\in\Theta_b$,
$\epsilon_n\equiv\dot{\ell}_a(Y;\alpha_{n}+\int_0^1 X(t)\beta_{n}(t)\,dt)$
satisfies $E\{\epsilon^4_n\mid X\}\le M_4$, a.s.
Then for any $\varepsilon>0$,
there exist positive constants $N_\varepsilon$ and $C_\varepsilon$
s.t. when $n\ge N_\varepsilon$,
\[
\inf_{\theta_n\in\Theta_b\dvtx  \llVert\theta_n\rrVert\ge
C_\varepsilon\eta_n} P_{\theta_n}(\mbox{reject }H_0)
\ge1-\varepsilon,
\]
where $\eta_n\asymp\sqrt{(nh^{1/2})^{-1}+\lambda}$.
\end{Theorem}

The model assumption $\inf_{y\in\mathcal{Y},a\in\bbR
}(-\ddot{\ell}_a(y;a))>0$ trivially holds for Gaussian regression and
exponential family considered in Examples~\ref{exa1pss}
and~\ref{exa2exponentialfamily}. As for the logistic model with
$L^2$ bounded $X$,
this condition can be replaced by $\inf_{y\in\mathcal{Y},a\in
\mathcal{I}}(-\ddot{\ell}_a(y;a))
>0$ under which the same conclusion as in Theorem~\ref{mrt} holds,
where $\mathcal{I}$ is some bounded open interval including the range
of $\langle R_X,\theta_{n0}\rangle$
for every $\theta_{n0}\in\Theta_b$.

Theorem~\ref{mrt} states that the PLRT is able to detect any local
alternative with separation rate no faster than $\eta_n$. In
particular, the minimum separation rate, that is, $n^{-2k/(4k+1)}$, is
achieved when $h=h^{\ast\ast}$. Note that $h^{\ast\ast}$ satisfies
the rate conditions required by Theorem~\ref{mrt}. For example, when
$k=m+r+1$, $a=r+1$, $r=0$
(see the discussions below Proposition~\ref{rateconv}), we can verify
this fact for $m>(7+\sqrt{33})/8\approx
1.593$ by direct calculations. In the specific $\ell_2$ regression,
Corollary~4.6 of \cite{HMV13} proves that the above minimax rate, that
is, $n^{-2k/(4k+1)}$, is optimal but under the perfect alignment
condition. Therefore, we prove that the proposed PLRT can achieve the
minimax optimality under more general settings.

The likelihood ratio testing procedure developed in this section
requires prior knowledge on the smoothness of the true slope function
and covariance kernel function, which might not be available in
practice. This motivates us to propose two adaptive testing procedures
in the next section.

\section{Adaptive testing construction}\label{secadaptivePLRT}

In this section, we develop two adaptive testing procedures for
$H_0\dvtx \beta=\beta_0$ without knowing
$m$ and $r$, that is, the smoothness of the true slope function and
covariance kernel function. One works for Gaussian errors, and another
works for sub-Gaussian errors. The test statistics for both cases are
maximizers over a sequence of standardized PLRTs. We derive the null
limit distribution as an extreme value distribution using Stein's
method \cite{CCK13,Stein86}. Their minimax properties will also be
carefully studied. To the best of our knowledge, our adaptive testing
procedures are the first ones developed in the roughness regularization
framework, which forms an interesting complement to those based on FPCA
techniques \cite{HMV13,L13}.

In this section, we focus on the $\ell_2$ regression with two types of
error: Gaussian error (Section~\ref{testgau}) and sub-Gaussian error
(Section~\ref{testsubg}). For simplicity, we assume $\beta_0=0$,
$\alpha=0$, and the errors to be of unit standard deviations. In
addition, we assume that the covariate process $X(t)$ has zero mean and
is independent of the error term. We remark that it is possible to
extend our results in this section to the general loss functions, but
with extremely tedious technical arguments.

Our test statistic is built upon a modified estimator $\widetilde
{\beta}_{n,\lambda}$ that is constructed in the following three
steps. The first step is to find a sequence\vspace*{1pt} of empirical eigenfunctions
$\widehat{\varphi}_\nu$ that satisfy $\widehat{V}(\widehat{\varphi
}_\nu,\widehat{\varphi}_\mu)=\delta_{\nu\mu}$ for all $\nu,\mu
\ge1$, where $\widehat{V}(\beta,\widetilde{\beta})=\int_0^1\int
_0^1\widehat{C}(s,t)\beta(s)\widetilde{\beta}(t)\,ds\,dt$ and $\widehat
{C}(s,t)=n^{-1}\sum_{i=1}^nX_i(s)X_i(t)$. We offer two methods for
finding $\widehat{\varphi}_\nu$. The\vspace*{1pt} first method conducts a
spectral decomposition, $\widehat{C}(s,t)=\sum_{\nu=1}^\infty
\widehat{\zeta}_\nu\widehat{\psi}_\nu(s)\widehat{\psi}_\nu(t)$,
with some nonincreasing positive sequence $\widehat{\zeta}_\nu$
and orthonormal functions $\widehat{\psi}_\nu$ in the usual
$L^2$-norm. Construct $\widehat{\varphi}_\nu=\widehat{\psi}_\nu
/\sqrt{\widehat{\zeta}_\nu}$. This method is easy to implement, but
implicitly assumes the perfect alignment condition. Our second method
is more generally applicable, but requires more tedious implementation.
Specifically, we apply similar construction techniques as in
Section~S.5 (Supplementary Material \cite{SCFDA}) by using the sample versions of
$\widetilde
{K}$, $C$ and $T$ therein. In particular, we choose $m=1, 2$ such that
the true slope function is more possible to be covered.

The\vspace*{1pt} second step is to define a data-dependent parameter space. Note
that $H^m(\mathbb I)$ can be alternatively defined as $\{\sum_{\nu
=1}^\infty b_\nu\varphi_\nu\dvtx  \sum_{\nu=1}^\infty b_\nu^2\nu
^{2k}<\infty\}$, where $k$ depends on $m$ in an implicit manner. An
approximate parameter space is $\mathcal{B}_k= \{\sum_{\nu
=1}^\infty b_\nu\widehat{\varphi}_\nu\dvtx  \sum_{\nu=1}^\infty b_\nu
^2\nu^{2k}<\infty\}$.
The consistency of the sample eigenfunctions implies that $\mathcal{B}_k$
is a reasonable approximation of $H^m(\mathbb I)$; see \cite{HMW06}.
The data-dependent parameter space is thus defined as
\[
\mathcal{B}_{kn}\equiv\Biggl\{\sum_{\nu=1}^n
b_\nu\widehat{\varphi}_\nu\bigg|  b_1,
\ldots,b_n\in\bbR\Biggr\}.
\]
In $\mathcal{B}_{kn}$, we can actually use the first $K_n\rightarrow
\infty$ ($K_n\ll n$) eigenfunctions as the basis. However, such a
general choice would give rise to unnecessary tuning of $K_n$ in practice.

In the last step, we obtain the desirable estimator as $\widetilde
{\beta}_{n,\lambda}=\break \arg\sup_{\beta\in\mathcal{B}_{kn}}\ell
_{n,\lambda}(\beta)$, where
%
\begin{equation}
\label{gauplik1} \ell_{n,\lambda}(\beta)=-\frac{1}{n}\sum
_{i=1}^n \Biggl(Y_i-\sum
_{\nu=1}^{n} b_\nu\omega_{i\nu}
\Biggr)^2\bigg/2-(\lambda/2)\sum_{\nu=1}^{n}
b_\nu^2\nu^{2k},
\end{equation}
and $\omega_{i\nu}=\int_0^1X_i(t)\widehat{\varphi}_\nu(t)\,dt$ for
$i=1,\ldots,n$ and $\nu\ge1$. The smoothing parameter $\lambda$
depends on both $n$ and $k$, denoted as $\lambda_k$. In particular, we
choose $\lambda_k$ as $c_0^{2k}n^{-4k/(4k+1)}(\log\log
{n})^{2k/(4k+1)}$ for some constant $c_0>0$ irrelevant to $k$. As will
be seen in later theorems, this choice yields the minimax optimality of
the adaptive testing. Define $Y=(Y_1,\ldots,Y_n)^T$,
$b=(b_1,\ldots,b_n)^T$, $\Lambda_k=\operatorname{diag}(1^{2k},
2^{2k}, \ldots, n^{2k})$,
$\Omega_i=(\omega_{i1},\ldots,\omega_{in})$
and $\Omega=(\Omega_1^T,\ldots,\Omega_n^T)^T$.
Hence we can rewrite $-2n\ell_{n,\lambda}(\beta)$ as
\[
(Y-\Omega b)^T(Y-\Omega b)+n\lambda_k b^T
\Lambda_k b,
\]
whose minimizer [equivalently, the maximizer of $\ell_{n,\lambda
}(\beta)$] is $\widehat{b}_{n,k}=(\Omega^T\Omega+n\lambda_k\Lambda
_k)^{-1}\Omega^TY$.
Note that $\Omega^T\Omega=nI_n$ by the following analysis: for any
$\nu,\mu\ge1$,
\begin{eqnarray*}
\sum_{i=1}^n\omega_{i\nu}
\omega_{i\mu}&=&\sum_{i=1}^n\int
_0^1X_i(t)\widehat{
\varphi}_\nu(t)\,dt\int_0^1X_i(s)
\widehat{\varphi}_\mu(s)\,ds
\\
&=&\int_0^1\int_0^1
\sum_{i=1}^nX_i(s)X_i(t)
\widehat{\varphi}_\nu(s)\widehat{\varphi}_\mu(t)\,ds\,dt
\\
&=&n\int_0^1\int_0^1
\widehat{C}(s,t)\widehat{\varphi}_\nu(s)\widehat{
\varphi}_\mu(t)\,ds\,dt=n\delta_{\nu\mu}.
\end{eqnarray*}
Therefore, $\widehat{b}_{n,k}=(nI_n+n\lambda_k\Lambda_k)^{-1}\Omega
^TY$ and $\widetilde\beta_{n,\lambda}=(\widehat\varphi_1,\ldots,\widehat
\varphi_n)\widehat b_{n,k}$.

In the above analysis, we implicitly assume $k$ to be known. However,
the value of $k$ is usually unavailable in practice. To resolve this
issue, we will conduct our testing procedure over a sequence of integer
$k$, that is, $\{1, 2,\ldots, k_n\}$, as will be seen in the next two
subsections. The full adaptivity of testing procedure is achieved when
we allow $k_n\rightarrow\infty$ so that the unknown $k$ can
eventually be captured by this sequence.

\subsection{Gaussian error}\label{testgau}

In this subsection, we denote the PLRT as $\mathrm{PLRT}_k\equiv\ell
_{n,\lambda}(0)-\ell_{n,\lambda}(\widetilde{\beta}_{n,\lambda})$
due to its\vspace*{1pt} dependence on $k$. By plugging in the above form of
$\widetilde\beta_{n,\lambda}$, we obtain
%
\begin{equation}
\mathrm{PLRT}_k=-\frac{1}{2n}Y^T\Omega(nI_n+n
\lambda_k\Lambda_k)^{-1}\Omega
^TY.\label{gauplrt1}
\end{equation}
We next derive a standardized version of $\mathrm{PLRT}_k$ under $H_0$. Define
$d_\nu(k)=1/(1+\lambda_k\rho_\nu(k))$, where $\rho_\nu(k)=\nu
^{2k}$, for any $\nu,k\ge1$. Under $H_0$, we have $Y=\epsilon
=(\epsilon_1,\ldots,\epsilon_n)^T$, and thus $-2n \mathrm{PLRT}_k=\sum_{\nu
=1}^{n}d_{\nu}(k)\eta_\nu^2$ for $\eta_1,\ldots,\eta_n\stackrel
{\mathrm{i.i.d.}}{\sim}N(0,1)$ by straightforward calculation. Hence, we have $E\{
-2n \mathrm{PLRT}_k\}=\sum_{\nu=1}^{n} d_\nu(k)$
and $\operatorname{Var}(-2n \mathrm{PLRT}_{k})=2\sum_{\nu=1}^{n} d_\nu^2(k)$.
The standardized version of $\mathrm{PLRT}_k$ can be written as
%
\begin{equation}
\label{deftauk} \tau_k=\frac{-2n\cdot \mathrm{PLRT}_k-\sum_{\nu=1}^{n}d_\nu(k)}{
(2\sum_{\nu=1}^{n}d_\nu(k)^2)^{1/2}}.
\end{equation}
Inspired by Theorem~\ref{limitCT}, $\tau_k$ is presumably of
standard normal distribution for any particular $k$. However, $k$ is
often unavailable in practice. As discussed previously, we shall
construct the adaptive testing based on a sequence of $\tau_k$ as
follows: (i)~define $AT_n^*=\max_{1\leq k\le k_n}\tau_k$, and (ii)
standardize $AT_n^\ast$ as
\[
AT_n=B_n\bigl(AT_n^*-B_n\bigr),
\]
where $B_n$ satisfies $2\pi B_n^2\exp(B_n^2)=k_n^2$; see \cite
{Hall79}. By Cram\'{e}r \cite{C46}, $B_n=\break \sqrt{2\log{k_n}}-\frac
{1}{2}(\log\log{k_n}+\log{4\pi})/\sqrt{2\log{k_n}}+O(1/\log
{k_n})\asymp\sqrt{2\log{k_n}}$
as $n$ becomes sufficiently large.

%
\begin{Theorem}\label{adapasym}
Suppose $k_n\asymp(\log{n})^{d_0}$, for some constant $d_0\in(0,1/2)$.
Then for any $\bar{\alpha}\in(0,1)$, we have under $H_0\dvtx \beta=0$,
\[
P(AT_n\le c_{\bar{\alpha}})\rightarrow1-\bar{\alpha}\qquad \mbox{as $n
\rightarrow\infty$},
\]
where $c_{\bar{\alpha}}=-\log(-\log(1-\bar{\alpha}))$.
\end{Theorem}

The proof of Theorem~\ref{adapasym} is mainly based on
Stein's leave-one-out method
\cite{Stein86} since under $H_0$, $\tau_k$ can be written as a sum of
independent random variables, that is, $\tau_k=\sum_{\nu=1}^n[d_\nu
(k)/s_{n, k}](\eta_\nu^2-1)$, where $s_{n,k}^2=2\sum_{\nu=1}^n
d_\nu(k)^2$.\vspace*{1pt}

In the end, we investigate the optimality of the proposed adaptive
testing procedure.
Consider the local alternative $H_{1n}\dvtx  \beta\in\mathcal{B}_{k,1}$, where
\[
\mathcal{B}_{k,1}\equiv\Biggl\{\sum_{\nu=1}^\infty
b_\nu\widehat{\varphi}_\nu\dvtx  \sum
_{\nu=1}^\infty b_\nu^2
\nu^{2k}\leq1 \Biggr\},
\]
for some fixed but \emph{unknown} integer $k\ge1$. For any real
sequence $\mathfrak{b}=\{b_\nu\}$ satisfying $\sum_{\nu=1}^\infty
b_\nu^2\nu^{2k}\leq1$, let $\beta_{\mathfrak{b}}=\sum_{\nu
=1}^\infty b_\nu\widehat{\varphi}_\nu$ be the alternative function
value, and let $P_{\mathfrak{b}}$ be the corresponding probability
measure. The following result shows that the adaptive test $AT_n$
achieves the optimal minimax rate (up to an logarithmic order), that
is, $\delta(n,k)\equiv n^{-2k/(4k+1)}(\log\log{n})^{k/(4k+1)}$, for
testing the hypothesis $H_0\dvtx \beta=0$,
with the alternative set being certain Sobolev ellipsoid $\mathcal
{B}_{k,1}$; see \cite{HMV13}.

Define $\llVert\mathfrak{b}\rrVert_{\ell^2}^2=\sum_{\nu
=1}^\infty b_\nu^2$ and $\llVert\mathfrak{b}\rrVert
_{k,\ell^2}^2=\sum_{\nu=1}^\infty b_\nu^2\rho_\nu(k)$.

%
\begin{Theorem}\label{minimaxadapt}
Suppose $k_n\asymp(\log{n})^{d_0}$, for some constant $d_0\in(0,1/2)$.
Then, for any $\varepsilon\in(0,1)$,
there exist positive constants $N_\varepsilon$ and $C_\varepsilon$
such that for any $n\ge N_\varepsilon$,
\[
\mathop{\inf_{\llVert\mathfrak{b}\rrVert_{\ell^2}\ge
C_\varepsilon\delta(n,k)}}_{\llVert\mathfrak{b}\rrVert
_{k,\ell^2}\le1}P_{\mathfrak{b}}(
\mbox{reject }H_0)\ge1-\varepsilon.
\]
\end{Theorem}

In Gaussian white noise models, Fan \cite{F96} and Fan and Lin \cite{FL98}
proposed an adaptive Neyman test based on multiple standardized test,
and derived the null limit distribution using the Darling--Erd\H{o}s
theorem. Theorems~\ref{adapasym}
and~\ref{minimaxadapt} can be viewed as extensions of such results to
functional data under Gaussian errors.\vadjust{\goodbreak} However, the
Darling--Erd\H{o}s theorem is no longer applicable in our setup due to
the difference in modeling and test construction. Instead, we employ
the Stein leave-one-out method. More interestingly, Stein's method can
even be applied to handle sub-Gaussian errors, as will be seen in
Section~\ref{testsubg}.

\subsection{Sub-Gaussian error}\label{testsubg}

In this subsection, we consider models with sub-Gaussian errors; that
is, there exists a positive constant $C_\epsilon$
such that $E\{\exp(t\epsilon)\}\le\exp(C_\epsilon t^2)$ for all
$t\in\bbR$. Further relaxation to the error term with finite fourth
moment is straightforward, but requires more stringent conditions on
the design.
For simplicity, we assume deterministic design, and suppose that
$X_i$'s satisfy the following moment condition:
%
\begin{equation}
\label{4thmom} \max_{1\le\nu\le n}\sum_{i=1}^n
\omega_{i\nu}^4=o\bigl(n^{8/5}(\log
\log{n})^{-14/5}\bigr).
\end{equation}
Recall that $\omega_{i\nu}=\int_{0}^1X_i(t)\widehat\varphi_{\nu
}(t)\,dt$ and is nonrandom under deterministic design.
Condition (\ref{4thmom}) implies that for any $\nu=1,\ldots,n$, the
magnitudes of $\omega_{1\nu},\ldots,\omega_{n\nu}$ should be
comparable given the restriction that $\sum_{i=1}^n\omega_{i\nu
}^2=n$. It\vspace*{1pt} rules out the situation that the sequence $\omega_{i\nu}$
is spiked at $i=\nu$, that is, $\omega_{\nu\nu}^2=n$ and $\omega
_{i\nu}=0$, for any $i\neq\nu$. This special situation essentially
gives rise to $\Omega=\sqrt{n}I$ such that $\mathrm{PLRT}_k$ defined in (\ref
{gauplrt1}) can be written as a scaled sum of independent centered
squares of the errors. The leave-one-out method employed in Theorem
\ref{adapasym} can handle this special case.

We first standardize $\mathrm{PLRT}_k$. The non-Gaussian assumption yields a
substantially different design matrix. Hence,
the scale factor is chosen to be different from the one used in
Section~\ref{testgau},
as described below. The standardized version is defined as
\[
\tilde{\tau}_k=\frac{-2n\cdot \mathrm{PLRT}_k-\sum_{\nu=1}^{n}d_\nu(k)}{
(2\sum_{i\neq j}a_{ij}^2(k))^{1/2}},
\]
where
$a_{ij}(k)$ is the $(i,j)$th entry of $A_k\equiv n^{-1}\Omega
(I_n+\lambda_k\Lambda_k)^{-1}\Omega^T$ for $1\le i,j\le n$. Note
that the scale factor in $\tilde{\tau}_k$, that is,
the term $(2\sum_{i\neq j}a_{ij}(k)^2)^{1/2}$, differs from the one in
$\tau_k$.
Technically, this new scale factor will facilitate the
asymptotic theory developed later in this section.
Let $AT_n^*=\max_{1\le k\le k_n}\tilde{\tau}_k$,
and $AT_n=B_n(AT_n^*-B_n)$, where $B_n$ satisfies $2\pi B_n^2\exp
(B_n^2)=k_n^2$.

%
\begin{Theorem}\label{adapasymnongauss}
Suppose $k_n\asymp(\log{n})^{d_0}$, for some constant $d_0\in(0,1/2)$.
Furthermore, $\epsilon$ is sub-Gaussian, and (\ref{4thmom}) holds.
Then for any $\bar{\alpha}\in(0,1)$, we have under $H_0\dvtx \beta=0$,
\[
P(AT_n\le c_{\bar{\alpha}})\rightarrow1-\bar{\alpha}\qquad\mbox{as $n
\rightarrow\infty$},
\]
where $c_{\bar{\alpha}}=-\log(-\log(1-\bar{\alpha}))$.
\end{Theorem}

The proof of Theorem~\ref{adapasymnongauss} is mainly
based on Stein's exchangeable pair method; see \cite{Stein86}.\vadjust{\goodbreak}

We conclude this subsection by showing that the proposed adaptive test
can still achieve the optimal minimax rate (up to a logarithmic order)
specified in \cite{HMV13}, that is, $\delta(n,k)$, even under
non-Gaussian errors. Recall that $\delta(n,k)$, $\llVert\mathfrak
{b}\rrVert_{\ell^2}$, $\llVert\mathfrak{b}\rrVert
_{k,\ell^2}$ and
$P_{\mathfrak{b}}$ are defined in Section~\ref{testgau}.

%
\begin{Theorem}\label{minimaxadaptnongauss}
Suppose $k_n\asymp(\log{n})^{d_0}$, for some constant $d_0\in(0,1/2)$.
Furthermore, $\epsilon$ is sub-Gaussian, and (\ref{4thmom}) holds.
Then, for any $\varepsilon\in(0,1)$,
there exist positive constants $N_\varepsilon$ and $C_\varepsilon$
such that for any $n\ge N_\varepsilon$,
\[
\mathop{\inf_{\llVert\mathfrak{b}\rrVert_{\ell^2}\ge
C_\varepsilon\delta(n,k)}}_{
\llVert\mathfrak{b}\rrVert_{k,\ell^2}\le1} P_{\mathfrak
{b}}(
\mbox{reject }H_0)\ge1-\varepsilon.
\]
\end{Theorem}

\section{Simulation study}\label{secsimulation}

In this section, we investigate the numerical performance of the
proposed procedures for inference.
We consider four different simulation settings. The settings in
Sections~\ref{setting1}--\ref{setting3}
are exactly the same as those in Hilgert et al. \cite{HMV13} and Lei \cite{L13}
so that we can fairly compare our testing results with theirs. We focus
on models with Gaussian error and choose $m=2$, that is, cubic spline.
Confidence interval in Section~\ref{secasympCI}, penalized
likelihood ratio test in Section~\ref{secPLRT} and adaptive testing
procedure in Section~\ref{testgau} are examined.
The setting in Section~\ref{setting4} is about functional linear
logistic regression.
Size and power of the PLRT test are examined.

\subsection{Setting 1}\label{setting1}
\begin{figure}

\includegraphics{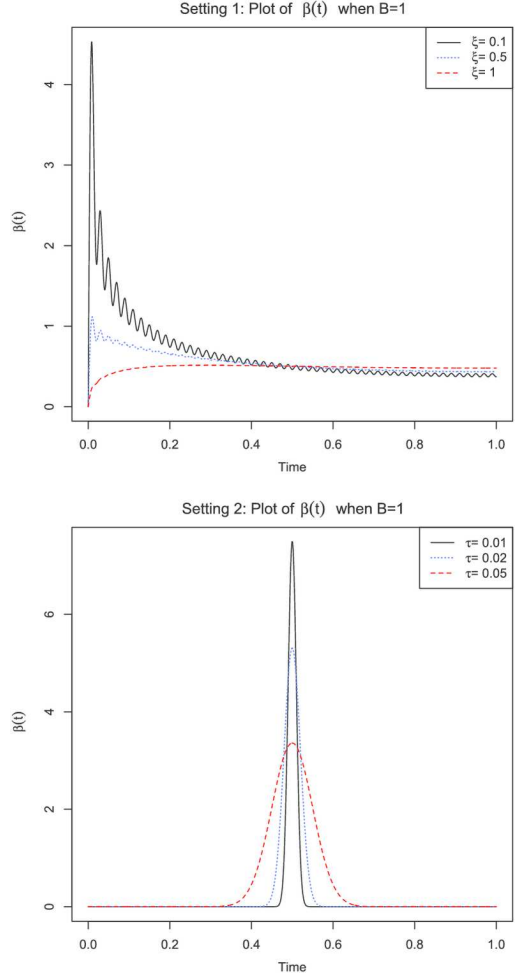}

\caption{Plots of $\beta_0(t)$ in settings 1 and 2.}\label{fig1}
\end{figure}

Data were generated in the same way as in Hilgert et al. \cite{HMV13}.
Consider the functional linear model $Y_i=\int_0^1 X_i(t)\beta_0(t)\,dt
+ \epsilon_i$, with $\epsilon_i$ being independent standard normal
for $i=1,\ldots,n$. Let $\lambda_j=(j-0.5)^{-2}\pi^{-2}$ and
$V_j(t)=\sqrt{2}\sin((j-0.5)\pi t)$, $t \in[0,1], j=1, 2, \ldots,
100$. The covariate curve $X_i(t)$ was Brownian motion simulated as
$X_i(t)=\sum_{j=1}^{100}\sqrt{\lambda_j}\eta_{ij}V_j(t)$, where
$\eta_{ij}$'s are independent standard normal for $i=1,\ldots,n$ and
$j=1,\ldots,100$. Each $X_i(t)$ was observed at $1000$ evenly spaced
points over $[0,1]$. The true slope function was chosen as
\[
\beta_0^{B,\xi}(t)= \frac{B}{\sqrt{\sum_{k=1}^{\infty}k^{-2\xi-1}}}\sum
_{j=1}^{100}j^{-\xi-0.5}V_j(t).
\]
Figure \ref{fig1} displays $\beta_0$. Four different signal strengths $B=(0,0.1,0.5,1)$ and three smoothness
parameters $\xi=(0.1,0.5,1)$
were considered. Note that $B=0$ implies $\beta_0=0$.

%

For each case study, we considered sample sizes $n=100$ and $n=500$
respectively,
and ran 10,000 trials to investigate the Monte Carlo performance of our methods.

\begin{longlist}
\item[\textit{Case study} 1: 95\% \textit{confidence interval for conditional mean}.]
In this study, we set $\mu_0(x_0) = E\{Y \mid X_0 = x_0\}=\int_0^1
x_0(t)\beta_0(t)\,dt$ with $B=1,\xi=1$,
where $x_0$ is independent of $X_1,\ldots,X_n$
and randomly generated from the same distribution as~$X_1$.
From (\ref{CIconditionalmean}),
the $95\%$ confidence interval for $\mu_0(x_0)$ is
\[
\bigl[\widehat Y_0 - n^{-1/2}z_{0.025}
\sigma_n,\widehat Y_0 + n^{-1/2}z_{0.025}
\sigma_n \bigr],
\]
where $\sigma_n^2=1+\sum_{\nu=1}^\infty x_\nu^2/(1+\lambda\rho
_\nu)$,
$x_\nu=\int_0^1 x_0(t)\varphi_\nu(t)\,dt$. Here $\varphi_\nu$ and
$\rho_\nu$ are both obtained through
(\ref{bvpexample}).

With 10,000 replications, percentages of the conditional mean $\mu_0(x_0)$
beyond the scope of CI and the average lengths of the confidence intervals
are summarized in Table~\ref{ci1}.

%
\begin{table}
\tabcolsep=0pt
\tablewidth=180pt
\caption{Case study 1: Percentage of $\mu_0(x_0)$ outside the $95 \%$
confidence intervals $\pm$ standard deviation
(average length of the $95 \%$ confidence intervals)}\label{ci1}
\begin{tabular*}{\tablewidth}{@{\extracolsep{\fill}}@{}lc@{}}
\hline
$\bolds{n=100}$ & $\bolds{n=500}$ \\
\hline
$4.89\pm0.42~(0.56)$ & $5.01\pm0.19~(0.39)$\\
\hline
\end{tabular*}
\end{table}

\item[\textit{Case study} 2: \textit{Size of the tests}.]
Denote\vspace*{1pt} the testing methods proposed by Hilgert et al. \cite
{HMV13} as $\mathrm{HMV13}^{(1)}$ and $\mathrm{HMV13}^{(2)}$. Under
$H_0\dvtx \beta=0$, we calculated the sizes of PLRT and AT (adaptive
testing), that is,\vspace*{1pt} the percentages of rejecting $H_0$, and then
compared them with $\mathrm{HMV13}^{(1)}$ and $\mathrm{HMV13}^{(2)}$
(directly cited from \cite{HMV13}) in Table~\ref{size}. Numerically,
we found that AT converges to the Gumbel distribution very slowly,
which is a common phenomenon in the extreme value literature; see \cite
{F96,FL98}. Following an idea similar to \cite{FL98}, finite sample
distributions of AT based on one million replications were instead
used. Obviously, from Table~\ref{size}, the proposed PLRT and AT are
both valid test statistics achieving desirable sizes.

%
\begin{table}[b]
\tabcolsep=0pt
\tablewidth=180pt
\caption{Case study 2: Sizes of the tests}\label{size}
\begin{tabular*}{\tablewidth}{@{\extracolsep{\fill}}@{}lcc@{}@{}}
\hline
& $\bolds{n=100}$ & $\bolds{n=500}$ \\
\hline
$\mathrm{HMV13}^{(1)}$ & $ 3.47~(\pm0.36)$ & $2.61~(\pm0.14)$\\
$\mathrm{HMV13}^{(2)}$ & $4.97~(\pm0.43)$ & $5.26~(\pm0.20)$\\
AT& $5.13~(\pm0.43)$ & $5.04~(\pm0.19)$ \\
PLRT & $5.45~(\pm0.45)$ & $5.19~(\pm0.20)$\\
\hline
\end{tabular*}
\end{table}

\item[\textit{Case study} 3: \textit{Power comparison}.]
In this study, we generated $\beta_0$ under different signal strengths
$B=(0.1,0.5,1)$ and smoothing parameters $\xi=(0.1,0.5,1)$.
Tables~\ref{power1} and~\ref{power2} summarize the powers of four
different testing methods, that is, the percentages of
rejecting $H_0\dvtx \beta=0$ at $95\%$ significance level, under $n=100$
and $n=500$. From $n=100$ to $n=500$, the powers of all tests increase.
In particular, PLRT generally performs better than AT since PLRT
incorporates known information from the model,
that is, $r=0$ (smoothness of the covariance kernel) and $m=2$
(smoothness of the functional parameter), while AT is adaptive on these
quantities. The power loss is the price paid for adaptiveness. We also
note that for weaker signals $B=0.1$, the\vspace*{1pt} powers of PLRT and AT improve
those of $\mathrm{HMV13}^{(1)}$, $\mathrm{HMV13}^{(2)}$, while for
stronger signals $B=0.5,1$, the powers of all tests are comparable.
\end{longlist}

%
\begin{table}
\tabcolsep=0pt
\caption{Case study 3: $n=100$. Powers}\label{power1}
\begin{tabular*}{\tablewidth}{@{\extracolsep{\fill}}@{}lca{5.7}a{5.7}a{5.7}@{\hspace*{-2pt}}}
\hline
& \textbf{Test} & \multicolumn{1}{c}{$\bolds{B=0.1}$} & \multicolumn{1}{c}{$\bolds{B=0.5}$} & \multicolumn{1}{c@{}}{$\bolds{B=1}$}\\
\hline
 $\xi=0.1$  &  $\mathrm{HMV13}^{(1)}$  &  3.88;(\pm0.38)  &  21.41;(\pm0.8)  &  77.24;(\pm0.82)  \\
&  $\mathrm{HMV13}^{(2)}$  &  5.80;(\pm0.46)  &  26.38;(\pm0.86)  & 81.78;(\pm0.76)  \\
& AT &  6.12;(\pm0.47 )  &  30.77;(\pm0.90 )  &  81.56;(\pm0.76 ) \\
& PLRT &  21.27;(\pm0.80) &  42.34;(\pm0.97)  &  84.20;(\pm0.71)
\\[3pt]
 $\xi=0.5$ &  $\mathrm{HMV13}^{(1)}$  &   4.74;(\pm0.42)  &  46.47;(\pm0.98)  &  98.68;(\pm0.22)  \\
&  $\mathrm{HMV13}^{(2)}$  &   6.65;(\pm0.49)  &  52.79;(\pm0.98)  & 99.06;(\pm0.19)  \\
& AT &  8.28;(\pm0.54 )  &  71.08;(\pm0.89 )  &   99.86;(\pm0.07)  \\
& PLRT &   23.13;(\pm0.83) &  74.74;(\pm0.85)  &  99.70;(\pm0.11)
\\[3pt]
 $\xi=1$  &  $\mathrm{HMV13}^{(1)}$  &   4.8;(\pm0.42)  &  62.67;(\pm0.95)  &  99.75;(\pm0.10)  \\
&  $\mathrm{HMV13}^{(2)}$  &   7.07;(\pm0.5)  &  68.30;(\pm0.91)  & 99.84;(\pm0.08)  \\
& AT &  9.47;(\pm0.57)  &  83.20;(\pm0.73 )  &  99.98;(\pm0.03) \\
& PLRT &   23.95;(\pm0.84)  &  84.03;(\pm0.72)  &  99.98;(\pm0.03)  \\
\hline
\end{tabular*}
\end{table}
%

%
\begin{table}[b]
\caption{Case study 3: $n=500$. Powers}\label{power2}
\begin{tabular*}{\tablewidth}{@{\extracolsep{\fill}}@{}lca{5.7}a{5.7}a{3.4}@{}}
\hline
& \textbf{Test} & \multicolumn{1}{c}{$\bolds{B=0.1}$} & \multicolumn{1}{c}{$\bolds{B=0.5}$} & \multicolumn{1}{c@{}}{$\bolds{B=1}$}\\
\hline
$\xi=0.1$ & $\mathrm{HMV13}^{(1)}$ &  5.17;(\pm0.19) & 86.98;(\pm0.29) & 100;(\pm0) \\
& $\mathrm{HMV13}^{(2)}$ &  8.48;(\pm0.24) & 90.89;(\pm0.25) &100;(\pm0) \\
& AT & 9.57;(\pm0.26) &89.14;(\pm0.27 ) &100;(\pm0)\\
& PLRT & 20.00;(\pm0.35) & 88.19;(\pm0.28) & 100;(\pm0)
\\[3pt]
$\xi=0.5$ & $\mathrm{HMV13}^{(1)}$ &  8.81;(\pm0.25) & 99.85;(\pm0.03) & 100;(\pm0) \\
& $\mathrm{HMV13}^{(2)}$ &  13.07;(\pm0.30) & 99.88;(\pm0.03) &100;(\pm0) \\
& AT &20.20;(\pm0.35 ) &100;(\pm0) &100;(\pm0)\\
& PLRT &  29.47;(\pm0.40) & 99.90;(\pm0.03) & 100;(\pm0)
\\[3pt]
$\xi=1$ & $\mathrm{HMV13}^{(1)}$ &  11.38;(\pm0.28) & 99.99;(\pm0.01) & 100;(\pm0) \\
& $\mathrm{HMV13}^{(2)}$ &  16.13;(\pm0.32) & 100;(\pm0) & 100;(\pm0) \\
& AT &26.51;(\pm0.39) &100;(\pm0 ) &100;(\pm0)\\
& PLRT &  34.08;(\pm0.42) & 100;(\pm0) & 100;(\pm0) \\
\hline
\end{tabular*}
\end{table}

\subsection{Setting 2}\label{setting2}
Let the true slope function be
\[
\beta_0^{B, \tau}(t)=B \exp\biggl\{-\frac{(t-0.5)^2}{2\tau
^2} \biggr
\} \biggl\{\int_0^1\exp\biggl\{-
\frac{(x-0.5)^2}{\tau
^2} \biggr\}\,dx \biggr\}^{-1/2},
\]
where $B=(0.5,1,2)$ and $\tau=(0.01,0.02,0.05)$. The processes
$X_i(t)$ and the samples were generated in the same way
as in Setting~1.\vadjust{\goodbreak}

The powers in Setting~2 are summarized in Tables~\ref{power3} and~\ref
{power4}. We observe similar phenomena as in Setting~1, that under
weaker signals, say $\tau=0.01, B=0.5$, PLRT and AT demonstrate larger
powers, while the powers of all procedures become comparable under
stronger signals. Again, PLRT generally has larger powers than the
adaptive procedure AT. All the powers increase as sample size becomes larger.

%
\begin{table}
\tabcolsep=0pt
\caption{Setting~2: $n=100$. Powers}\label{power3}
\begin{tabular*}{\tablewidth}{@{\extracolsep{\fill}}@{}lca{5.7}a{5.7}a{5.7}@{\hspace*{-2pt}}}
\hline
& \textbf{Test} & \multicolumn{1}{c}{$\bolds{B=0.5}$} & \multicolumn{1}{c}{$\bolds{B=1}$} & \multicolumn{1}{c@{}}{$\bolds{B=2}$}\\
\hline
$\tau=0.01$ & $\mathrm{HMV13}^{(1)}$ & 4.94;(\pm0.42) & 11.85;(\pm0.63) & 46.69;(\pm0.98) \\
& $\mathrm{HMV13}^{(2)}$ & 7.25;(\pm0.51) & 15.49;(\pm0.71) &53.56;(\pm0.98) \\
& AT & 9.88;(\pm0.58) & 23.86;(\pm0.84 ) & 69.46;(\pm0.90 )\\
& PLRT & 17.9;(\pm0.75 )& 33.25;(\pm0.92 ) & 81.04;(\pm0.77)
\\[3pt]
$\tau=0.02$& $\mathrm{HMV13}^{(1)}$ &  7.33;(\pm0.51) & 23.09;(\pm0.83) & 80.26;(\pm0.78) \\
& $\mathrm{HMV13}^{(2)}$ &  10;(\pm0.59) & 28.54;(\pm0.89) &84.04;(\pm0.72) \\
& AT & 14.58;(\pm0.69) & 42.21;(\pm0.97) & 93.54;(\pm0.48) \\
& PLRT &  22.87;(\pm0.82 )& 53.21;(\pm0.98 ) & 97.83;(\pm0.29)
\\[3pt]
$\tau=0.05$ & $\mathrm{HMV13}^{(1)}$ &  13.85;(\pm0.68) & 56.51;(\pm0.97) & 99.48;(\pm0.14) \\
& $\mathrm{HMV13}^{(2)}$ &  18.13;(\pm0.50) & 63.09;(\pm0.95) &99.65;(\pm0.12) \\
& AT & 28.31;(\pm0.88) & 78.52;(\pm0.80 ) & 99.96;(\pm0.04)\\
& PLRT &  37.54;(\pm0.95) & 87.63;(\pm0.65) & 100;(\pm0)\\
\hline
\end{tabular*}\vspace*{-6pt}
\end{table}
%

%
\begin{table}[b]
\tabcolsep=0pt
\caption{Setting~2: $n=500$. Powers}\label{power4}
\begin{tabular*}{\tablewidth}{@{\extracolsep{\fill}}@{}lca{5.7}a{5.7}a{5.7}@{\hspace*{-2pt}}}
\hline
& \textbf{Test} & \multicolumn{1}{c}{$\bolds{B=0.5}$} & \multicolumn{1}{c}{$\bolds{B=1}$} & \multicolumn{1}{c@{}}{$\bolds{B=2}$}\\
\hline
$\tau=0.01$ & $\mathrm{HMV13}^{(1)}$ & 12.41;(\pm0.42) & 54.6;(\pm0.63) & 99.75;(\pm0.98) \\
& $\mathrm{HMV13}^{(2)}$ & 17.99;(\pm0.51) & 63.16;(\pm0.71) &99.98;(\pm0.98) \\
& AT & 28.93;(\pm0.40 ) & 79.75;(\pm0.35 ) & 100;(\pm0 )\\
& PLRT & 34.77;(\pm0.42 )&  86.08;(\pm0.30 ) & 100;(\pm0 )
\\[3pt]
$\tau=0.02$& $\mathrm{HMV13}^{(1)}$ &  26.11;(\pm0.51) &88.91;(\pm0.83) & 100;(\pm0) \\
& $\mathrm{HMV13}^{(2)}$ &  33.95;(\pm0.59) & 92.62;(\pm0.89) &100;(\pm0) \\
& AT & 50.25;(\pm0.44 ) & 97.03;(\pm0.15) & 100;(\pm0) \\
& PLRT &  56.57;(\pm0.43)& 99.20;(\pm0.08) & 100;(\pm0)
\\[3pt]
$\tau=0.05$ & $\mathrm{HMV13}^{(1)}$ &  65.38;(\pm0.68) & 99.95;(\pm0.97) & 100;(\pm0) \\
& $\mathrm{HMV13}^{(2)}$ &  72.74;(\pm0.50) & 99.99;(\pm0.95) &100;(\pm0) \\
& AT & 86.92;(\pm0.30) & 100;(\pm0 ) & 100;(\pm0)\\
& PLRT &  92.07;(\pm0.24) & 100;(\pm0 ) & 100;(\pm0 ) \\
\hline
\end{tabular*}\vspace*{-7pt}
\end{table}
%

\subsection{Setting~3}\label{setting3}
In this setting, data were generated in the same way as in Section~4.2
of \cite{L13}. Hence we will compare our PLRT and AT with the testing
procedure in \cite{L13}, denoted as L13. Specifically, the covariance
operator has eigenvalues $\kappa_j = j^{-1.7}$ and eigenfunctions
$\phi_1(t)=1, \phi_j(t)=\sqrt{2}\cos((j-1)\pi t)$ for $j \geq2$.
The covariate processes are $X_i(t)=\sum_{j=1}^{100}\sqrt{\kappa
_j}\eta_j \phi_j(t)$, where $\eta_j$'s are independent standard
normal. Each $X_i(t)$ was observed on $1000$ evenly spaced points over~$[0,1]$.

In the first case denoted as $\operatorname{Model}(2,1)$, let $\theta_j = \bar
{\theta}_j/\llVert\bar{\theta}\rrVert_2$,
where $\bar{\theta}_j=0$ for $j>2$, $\bar{\theta}_j=b_j\cdot I_j$
for $j=1,2$, $b_1$ and $b_2$
are independent $\operatorname{Unif}(0,1)$, and\vadjust{\goodbreak} $(I_1,I_2)$ follows a multinomial
distribution $\operatorname{Mult}(1;0.5,0.5)$.
Let the true function be
$\beta_0(t)=r\sum_{j=1}^{100}\theta_j\phi_j(t)$, where
$r^2=(0,1,0.2,0.5,1.5)$.

%
\begin{table}
\tabcolsep=0pt
\caption{Setting~3: Powers}\label{power5}
\begin{tabular*}{\tablewidth}{@{\extracolsep{\fill}}@{}lcca{5.7}a{5.7}a{5.7}a{5.7}@{\hspace*{-2pt}}}
\hline
& \textbf{Sample size} & \textbf{Test} &
\multicolumn{1}{c}{$\bolds{r^2=0.1}$} &
\multicolumn{1}{c}{$\bolds{r^2=0.2}$} &
\multicolumn{1}{c}{$\bolds{r^2=0.5}$} &
\multicolumn{1}{c@{}}{$\bolds{r^2=1.5}$}\\
\hline
$\operatorname{Model} (2,1)$ & $n=50$ & L13 & 16.20;(\pm1.02)  &  26.40;(\pm1.22) & 54.20;(\pm1.38) & 80.80;(\pm1.09) \\
& & AT & 47.81;(\pm1.38) &  64.94;(\pm1.32 )&  84.75;(\pm1.00) & 99.13;(\pm0.26 )\\
& & PLRT & 57.76;(\pm1.37) &  72.70;(\pm1.23) & 90.18;(\pm0.82) & 99.52;(\pm0.19) \\[3pt]
& $n=100$ & L13 & 25.80;(\pm0.86 ) &  42.20;(\pm0.97) & 68.20;(\pm0.91) & 90.40;(\pm0.58) \\
& & AT & 65.53;(\pm0.93 ) &  79.75;(\pm0.79 )&  96.97;(\pm0.34) & 99.99;(\pm0.02 ) \\
& & PLRT & 74.04;(\pm0.86 ) &  87.98;(\pm0.64) & 98.22;(\pm0.26) & 100;(\pm0) \\[3pt]
& $n=500$ & L13 & 67.20;(\pm0.41) &  84.60;(\pm0.32) & 94.40;(\pm0.20) & 97.20;(\pm0.14) \\
& & AT & 97.81;(\pm0.13 ) &  100;(\pm0)&  100;(\pm0 ) & 100;(\pm0 ) \\
& & PLRT & 98.5;(\pm0.11) &  99.94;(\pm0.02) & 100;(\pm0) &100;(\pm0)
\\[6pt]
$\operatorname{Model} (9,2)$ & $n=50$ & L13 & 9.00;(\pm0.79) &  14.00;(\pm0.96)& 29.60;(\pm1.27) & 43.40;(\pm1.37) \\
& & AT & 21.72;(\pm1.14) &  27.57;(\pm1.24 )&  37.67;(\pm1.34) & 53.33;(\pm1.38 )\\
& & PLRT & 39.54;(\pm1.36 ) &  46.22;(\pm1.38) & 56.92;(\pm1.37) & 73.42;(\pm1.22) \\[3pt]
& $n=100$ & L13 & 13.4;(\pm0.67) &  27.8;(\pm0.88) & 39.8;(\pm0.96) & 65.8;(\pm0.93) \\
& & AT & 27.86;(\pm0.88 ) &  21.63;(\pm0.81)&  47.61;(\pm0.98) & 65.69;(\pm0.93) \\
& & PLRT & 45.80;(\pm0.98) &  53.72;(\pm0.98) & 67.12;(\pm0.92) & 83.88;(\pm0.72) \\[3pt]
& $n=500$ & L13 & 42.40;(\pm0.43) &  47.8;(\pm0.44) & 72.4;(\pm0.39) & 93.4;(\pm0.22) \\
& & AT & 49.40;(\pm0.44 ) &  58.29;(\pm0.43)&  80.21;(\pm0.35) & 91.23;(\pm0.25 ) \\
& & PLRT & 69.44;(\pm0.40) &  80.00;(\pm0.35) & 91.70;(\pm0.24) & 99.22;(\pm0.08) \\
\hline
\end{tabular*}
\end{table}

In the second case denoted as $\operatorname{Model}(9,2)$, a different choice of
$\theta_j$ was considered.
Specifically, $\theta_j = \bar{\theta}_j/\llVert\bar{\theta
}\rrVert_2$, where $\bar{\theta}_j=0$ for $j>9$,
$\bar{\theta}_j=b_j\cdot I_j$ for $j=1,\ldots,9$, $b_1,\ldots,b_9$
are independent $\operatorname{Unif}(0,1)$, and $(I_1,\ldots,I_9)$ follows a
multinomial distribution $\operatorname{Mult}(2;1/9,\ldots,1/9)$.

In both cases, the samples were drawn from $Y_i=\int_{0}^1 X_i(t)\beta
(t)\,dt + \epsilon_i$, $i=1,\ldots,n$,
where $\epsilon_i$ are independent standard Gaussian. $5000$ Monte
Carlo trials
were conducted in each case under different sample sizes $n=50, 100$
and $500$.

Results are summarized in Table~\ref{power5},
from which we can see that the powers of AT and PLRT improve those
of L13, especially when $r^2=0.1,0.2$ (weaker signals). As $n$ increases,
the power of L13 becomes more comparable to those of PLRT and AT
especially when $r^2=1.5$ (stronger signal).
Again, PLRT generally has larger powers than adaptive methods.

\subsection{Setting~4}\label{setting4}

Let $Y\in\{0,1\}$ be a binary variable
generated from the following functional logistic regression model:
\[
P(Y=1\mid X)=\frac{\exp(\int_0^1X(t)\beta_0(t)\,dt)}{1+\exp(\int
_0^1X(t)\beta_0(t)\,dt)}.
\]
The predictor process $X_i$ was simulated as $X_i(t)=\sum
_{j=1}^{100}\sqrt{\lambda_j}\eta_{ij}V_j(t)$,
where $\lambda_j$ and $V_j(t)$ are exactly the same as in Setting~1,
$\eta_{ij}$'s are independent truncated normals,
that is, $\eta_{ij}=\xi_{ij}I_{\{\llvert \xi_{ij}\rrvert \leq0.5\}} +
0.5I_{\{\xi
_{ij}> 0.5\}} -0.5I_{\{\xi_{ij}< -0.5\}}$,
with $\xi_{ij}$ being a standard normal random variable.
Each $X_i(t)$ was observed at $1000$ evenly spaced points over $[0,1]$.
We intend to test $H_0\dvtx  \beta=0$.
To examine the power, data were generated under $\beta
_0(t)=3*10^5(t^{11}(1-t)^6)$
for $t\in[0,1]$.

We examined two sample sizes: $n=100$ and $n=500$.
Results (summarized in Table~\ref{setting4sizepower})
were based on 10,000 independent trials.
It can be seen that when $n=100$ and $500$,
the test achieves the desired sizes.
The power at $n=100$ is small,
but the power at $n=500$ approaches one, demonstrating the asymptotic
property of the test.

%
\begin{table}
\tabcolsep=0pt
\tablewidth=180pt
\caption{Setting~4: Size and power}\label{setting4sizepower}
\begin{tabular*}{\tablewidth}{@{\extracolsep{\fill}}@{}lcc@{}}
\hline
& $\bolds{n=100}$ & $\bolds{n=500}$ \\
\hline
Size & 0.054 & 0.046 \\
Power& 0.387 & 0.985 \\
\hline
\end{tabular*}
\end{table}

\section{Discussion}\label{secdisc}

The current paper and our previous work on nonparametric regression
models \cite{SC13} are both built upon the RKHS framework and theory.
Hence it seems necessary for us to comment their technical connections
and differences to facilitate the reading. Compared to \cite{SC13},
the RKHS considered in the current paper has a substantially different
structure that involves a covariance function of the predictor process.
This immediately causes a difference in building the eigensystems:
\cite{SC13} relies on an ODE system, but the current paper relies on
an integro-differential system. Hence the methods of analyzing both
systems are crucially different. Meanwhile, the asymptotic analysis on
the statistical inference such as the penalized likelihood ratio test
are also different. For example, \cite{SC13} only considers the
reproducing kernel, while the current work requires a delicate
interaction between
the reproducing kernel and the covariance kernel. More importantly, the
relaxation of perfect alignment between both kernels poses more
technical challenges.

Besides, Assumption~\ref{A3} requires $\llVert\varphi_\nu\rrVert
_{L^2}\le C_\varphi\nu^a$ for $\nu\ge1$ and a constant $a\ge
0$. The introduction of factor $a$ in Assumption~\ref{A3} is helpful
in simplifying our proofs. However,
it is interesting to investigate how to avoid imposing this seemingly
``redundant'' $a$.
As indicated by Proposition~\ref{validA3}, that $a$ relates to $C$
(and hence $V$), one possible strategy is to avoid the use of $V$.
Instead, one may use its empirical version, namely $V_n$, as suggested
by one referee. This would require a delicate analysis of the
convergence of $V_n$,
which may be handled by techniques in \cite{M10}. We leave this as a
future exploration.


\section*{Acknowledgments}
We thank Pang Du for providing us his R
code for functional logistic regression,
and thank Ph.D. student Meimei Liu at Purdue for help with the simulation
study. We also thank Co-editor
Runze Li, an Associate Editor, and two referees for helpful comments
that lead to
important improvements on the paper.

Guang Cheng was on sabbatical at Princeton while the revision of this
work was carried out; he would like
to thank the Princeton ORFE department for its hospitality and support.

\begin{supplement}[id=suppA]
\stitle{Supplement to ``Nonparametric inference in generalized functional~linear models''}
\slink[doi]{10.1214/15-AOS1322SUPP} 
\sdatatype{.pdf}
\sfilename{aos1322\_supp.pdf}
\sdescription{Proofs are provided.}
\end{supplement}

%

\printaddresses
\end{document}